\documentclass[a4paper, 12pt]{amsart}
\usepackage{amssymb,amsmath,latexsym,amsthm}
\usepackage[latin1]{inputenc}
\usepackage[T1]{fontenc}
\usepackage[english]{babel}
\usepackage{amssymb}
\usepackage{amsmath}
\usepackage{amsthm}
\usepackage{amscd}
\usepackage{amsfonts}
\usepackage{stmaryrd}
\usepackage{pb-diagram}
\usepackage{epic,eepic,epsfig}
\usepackage{a4wide}
\usepackage{xypic}
\usepackage{verbatim}
\usepackage{nextpage}
\usepackage{fancyhdr}
\usepackage[OT2,T1]{fontenc}
\DeclareSymbolFont{cyrletters}{OT2}{wncyr}{m}{n}
\DeclareMathSymbol{\Sha}{\mathalpha}{cyrletters}{"58}
\long\def\proof#1{\removelastskip\vskip\baselineskip\relax\noindent{\it
Proof\if!#1!\else\ \ignorespaces#1\fi.\ }\ignorespaces}
\renewcommand{\i}{\iota}

\newcommand{\wh}{\widehat}

\newcommand{\lgc}[2]{\mbox{$\bigl(\frac{#1}{#2}\bigr)_3$}}

\newcommand{\isom}{\simeq}
\newcommand{\LR}{\longrightarrow}

\DeclareMathOperator{\Ker}{Ker}
\newcommand{\Ima}{\mbox{\rm Im}}

\DeclareMathOperator{\Gal}{Gal}
\DeclareMathOperator{\Sel}{Sel}
\DeclareMathOperator{\rk}{rk}

\DeclareMathOperator{\Aut}{Aut}

\newcommand{\ov}[1]{\overline{#1}}
\newcommand{\Q}{{\mathbb Q}}
\newcommand{\Z}{{\mathbb Z}}
\newcommand{\R}{{\mathbb R}}
\newcommand{\F}{{\mathbb F}}

\renewcommand{\v}{{\mathfrak v}}

\newcommand{\bv}{{\mathbf v}}

\newcommand{\la}{\lambda}

\renewcommand{\a}{{\mathfrak a}}

\newcommand{\f}{{\mathfrak f}}

\newcommand{\al}{\alpha}
\newcommand{\be}{\beta}
\newcommand{\ga}{\gamma}

\newcommand{\p}{{\mathfrak p}}
\newcommand{\q}{{\mathfrak q}}

\newcommand{\gd}{{\mathfrak d}}

\newcommand{\Proof}{{\it Proof. \/}}
\newcommand{\squareforqed}{\hbox{\rlap{$\sqcap$}$\sqcup$}}
\renewcommand{\qed}{\ifmmode\squareforqed\else{\unskip\nobreak\hfil
\penalty50\hskip1em\null\nobreak\hfil\squareforqed
\parfillskip=0pt\finalhyphendemerits=0\endgraf}\fi}

\newcommand{\fp}{\qed\removelastskip\vskip\baselineskip\relax}

\newtheorem{theorem}{Theorem}[section]
\newtheorem{corollary}[theorem]{Corollary}
\newtheorem{proposition}[theorem]{Proposition}
\newtheorem{lemma}[theorem]{Lemma}
\newtheorem{definition}[theorem]{Definition}

\newcommand{\litem}{\par\noindent\dimen0=\parindent%
    \advance\dimen0 by-4pt
               \hangindent=\dimen0\ltextindent}

\newcommand{\ltextindent}[1]{\hbox to \hangindent{#1\hss}\ignorespaces}
\newcommand{\ltextjndent}[1]{\hbox to \hangindent{#1\hss}\ignorespaces\kern-1ex}

\renewcommand{\pmod}[1]{\allowbreak\ ({\rm{mod}}\,\,#1)}

\begin{document}
\baselineskip=17pt

\def\refname{\centerline{Bibliography}}

\title{Elementary $3$-Descent with a $3$-Isogeny}

\author[Henri {\sc Cohen}]{{\sc Henri} Cohen}
\address{Henri {\sc Cohen}\\
IMB Universit\'e Bordeaux I\\
351 Cours de la Lib\'eration\\
33405 Talence, France}
\email{cohen@math.u-bordeaux1.fr}
\urladdr{http://www.math.u-bordeaux.fr/~cohen/}

\author[Fabien {\sc Pazuki}]{{\sc Fabien} Pazuki}
\address{Fabien {\sc Pazuki}\\
IMJ Universit\'e Paris 7\\
site Chevaleret\\
2, place de Jussieu\\
75 251 Paris Cedex 05, France}
\email{pazuki@math.jussieu.fr}
\urladdr{http://www.math.jussieu.fr/~pazuki}

\maketitle

\begin{abstract}
In this expository paper, we show how to use in practice $3$-descent with
a $3$-isogeny to find an estimate for the rank of an elliptic curve having
a rational $3$-torsion subgroup, and we also give a geometric interpretation
of these computations.
\end{abstract}

{\flushleft
\textbf{Keywords :} Elliptic curves, $3$-descent.\\
\textbf{Mathematics Subject Classification :} 11G05, 14G05, 14H52.}

\section{Introduction}

The aim of this work is to give a very explicit way to estimate the rank of an elliptic curve over $\mathbb{Q}$ using $3$-descent. We will suppose that the elliptic curve has a rational $3$-torsion subgroup. It allows us to pick an affine model of the form $y^2=x^3+D(ax+b)^2$. After introducing the descent maps, we explain in section 2 how to us $3$-descent. Then we show how to compute principal homogeneous spaces in the case $D=1$ in section 3. We do the same in section 4 for the case $D\neq 1$, which is a bit more technical. Sections 5 and 6 include all the results needed for local solubility. For the sake of brevity we do not include all the details of the calculations, they are of course availaible upon request. In section 7, one finds several examples of families where we find the $\mathbb{Q}$-rank, and we also give some applications, such as prime values of certain cubic forms.

The main strategy here is improving on \cite{Cohen3}, Section 8.4, where this explicit way of doing descent is explained for $2$-descent. The $3$-Selmer group has also been studied in \cite{Top}, precised by \cite{DeLo}. See for example \cite{CreFiOSiSto1, CreFiOSiSto2, Fish} and their references for a more general treatment.

We would like to thank J. Cremona for his interest and the anonymous referee for his advice.

\subsection{The Geometric Setting}

We recall a few facts about descent on elliptic curves. Let $E/k$ an elliptic
curve over a number field $k$ and let $n\geq 2$ be an integer. First, using
Galois cohomology, we have the short exact sequences:

\[
\xymatrix{
0 \ar[r] &\displaystyle{ E(k)/nE(k) \ar[d]\ar[r]^{\delta} }& H^{1}(k,E[n]) \ar[d]^{\prod_{v}\textrm{res}_{v}}\ar[r]\ar[rd]^{\varphi} &\displaystyle{ H^{1}(k,E)[n] \ar[d]^{\prod_{v}\textrm{res}_{v}}\ar[r] }& 0\\
0 \ar[r] &\displaystyle{ \prod_{v}E(k_{v})/nE(k_{v}) \ar[r]^{\prod_{v}\delta_{v}} }& \displaystyle{\prod_{v}H^{1}(k_{v},E[n]) \ar[r] }& \displaystyle{\prod_{v}H^{1}(k_{v},E)[n] \ar[r] }& 0\\}
\]

We recall the definition of the \emph{n-Selmer group}:
$$\textrm{Sel}^{\textrm{(n)}}(k,E):=\Ker\Big( \varphi : H^{1}(k,E[n])\LR \prod_{v}H^{1}(k_{v},E)[n]\Big).$$
We also recall the definition of the \emph{Tate--Shafarevich group}:
$$\Sha(k,E):=\Ker\Big(H^{1}(k,E)\LR \prod_{v}H^{1}(k_{v},E)\Big).$$
This leads to the following short exact sequence:
$$0\LR E(k)/nE(k) \LR \textrm{Sel}^{\textrm{(n)}}(k,E) \LR \Sha(k,E)[n] \LR 0 \;,$$
\noindent
where one can show that every term is a finite group, so that
$$\Big|\textrm{Sel}^{\textrm{(n)}}(k,E)\Big|=\Big|E(k)/nE(k)\Big|\,\Big|\Sha(k,E)[n]\Big|\;,$$
which gives
$$n^{\displaystyle{\rk(E/k)}}=\frac{\displaystyle{\Big|\textrm{Sel}^{\textrm{(n)}}(k,E)\Big|}}{\displaystyle{\Big|E(k)_{\textrm{tors}}/nE(k)_{\textrm{tors}}\Big|\,\Big|\Sha(k,E)[n]\Big|}}\;.$$
Thus, to get the exact value of the rank $\rk(E/k)$, we must compute the 
$n$-Selmer group and the $n$-torsion part of the Tate--Shafarevich group. 

\vspace{0.3cm}

Recall that a \emph{twist} of an object $X$ defined over $k$ is an object $Y$ 
defined over $k$ that is isomorphic to $X$ over $\bar{k}$.

Since $\textrm{Sel}^{\textrm{(n)}}(k,E)\subset H^{1}(k,E[n])$, if we find a 
geometric object $X$ such that $\Aut_{\bar{k}}(X)\cong E[n]$, we can 
interpret the elements of the $n$-Selmer group as twists of the object $X$.
This idea gives rise to different interpretations of the elements of the 
$n$-Selmer group, as is clearly explained in \cite{CreFiOSiSto1, CreFiOSiSto2}.
In the present paper we are going to describe explicitly the geometrical 
interpretation of those elements that we now recall. First, if $O$ denotes
the identity element of $E$, the complete linear system given by $|n.O|$ 
induces a morphism: $E\LR \mathbb{P}^{n-1}.$

\begin{definition} A \emph{diagram} $[C\LR S]$ is a morphism from a torsor $C$
under $E$ to a variety $S$. We will say that two diagrams $[C_{1}\LR S_{1}]$
and $[C_{2}\LR S_{2}]$ are isomorphic if the following diagram is commutative:

\[
\xymatrix{
C_{1}\ar[r]\ar[d]^{\cong} & S_{1} \ar[d]^{\cong}\\
C_{2} \ar[r] & S_{2}\\}
\]

We will define a \emph{Brauer--Severi diagram} $[C\LR S]$ to be a twist of the
diagram $X=[E\LR \mathbb{P}^{n-1}]$. In particular $S$ is a twist of 
$\mathbb{P}^{n-1}$, called a Brauer--Severi variety.
\end{definition}

\vspace{0.1cm}

Following \cite{CreFiOSiSto2}, we interpret an element of the $n$-Selmer group
of $E$ as a Brauer--Severi diagram $[C\LR S]$ such that the curve $C$ has 
points everywhere locally, hence one can take $S=\mathbb{P}^{n-1}$.

\vspace{0.3cm}

We now specialize to the case $n=3$. In this particular case, the Brauer--Severi 
diagrams we are looking for are of the type $[C\LR \mathbb{P}^{2}]$, the curve
$C$ being a plane cubic with points everywhere locally, given with an action 
of $E[3]$ on it by linear automorphisms.

\subsection{The Arithmetic Setting}

For the proofs of all the results given in this section, we refer to
\cite{Cohen3}, Section 8.4, although there are other pointers in the
literature. The $3$-Selmer group in this particular case has also been studied
in \cite{Top}, precised by \cite{DeLo}.

Let $E$ be an elliptic curve defined over $\Q$ and having a rational
$3$-torsion subgroup that we denote $\{O,T,-T\}$. Let us stress here that the point $T$ does not need to be rational itself. It is easy to see that $E$ can be given by
an affine equation of the type
$$y^2=x^3+D(ax+b)^2$$
with $a$, $b$, and $D$ in $\Q$, and the discriminant of $E$ is equal to
$16b^3D^2(4Da^3-27b)$, so we must have $b$ and $D$ nonzero and $4Da^3-27b\ne0$.
The $3$-torsion point $T$ is equal to $(0,b\sqrt{D})$, so is rational if and
only if $D\in{\Q^*}^2$.

\begin{lemma}\label{lemeq} There exists a unique equation of $E$ of the form
$y^2=x^3+D(ax+b)^2$, where $a$, $b$, and $D$ are in $\Z$, $D$ is a fundamental
discriminant (including $1$), $b>0$, and if we write $b=b_1b_3^3$ with $b_1$
cubefree then $(a,b_3)=1$.\end{lemma}

\Proof We can write uniquely $D=D_0f^2$, where $D_0$ is a fundamental
discriminant and $f\in\Q^*$, so our initial equation can be written
$y^2=x^3+D_0(a'x+b')^2$ with $a'=fa$ and $b'=fb$. Changing $(x,y)$ into 
$(x/u^2,y/u^3)$ changes $(a',b')$ into $(ua',u^3b')$, so it is clear that 
we may assume that $a$ and $b$ are in $\Z$, and if $b=b_1b_3^3$ and 
$g=(b_3,a)$, changing $(x,y)$ into $(xg^2,yg^3)$ changes $(a,b)$ into 
$(a/g,b/g^3)$, hence $(a,b_3)$ into $(a/g,b_3/g)$, so insures that
$(a,b_3)=1$. Finally, changing $(a,b)$ into 
$(-a,-b)$ insures that $b>0$. Uniqueness is immediate and left to the reader.
\fp

From now on, we will always assume that the equation of our curve is given
satisfying the conditions of the above lemma (although we will mainly use the
fact that $D$ is a fundamental discriminant), and we will denote by $K$
the field $K=\Q(\sqrt{D})$ of discriminant $D$, which will be equal to $\Q$
if $D=1$ and to a quadratic field otherwise.

Given such a curve $E$, our aim is to give an estimate for the rank of $E$,
and if possible the rank itself, using $3$-descent. We first recall the
definition and main properties of the $3$-descent maps.

\begin{definition} Let $E$ be an elliptic curve defined over $\Q$, choose an affine model given by an 
equation $y^2=x^3+D(ax+b)^2$ with $D$ a fundamental discriminant, and
let $T=(0,b\sqrt{D})$ be a $3$-torsion point.
\begin{enumerate}\item The $3$-descent map $\al$ is a map from
$E(\Q)$ to the subgroup $G_3$ of classes of elements of $K^*/{K^*}^3$
whose norm is a cube (or $G_3=\Q^*/{\Q^*}^3$ if $D=1$) defined by
$\al(O)=1$, $\al((0,b))=1/(2b)$ when $D=1$, and in general by
$\al((x,y))=y-(ax+b)\sqrt{D}$.
\item The curve $\wh{E}$ is defined by a similar equation 
$y^2=x^3+\wh{D}(\wh{a}x+\wh{b})^2$, where $\wh{D}=-3D$, $\wh{a}=a$, and
$\wh{b}=(27b-4a^3D)/9$, and the corresponding $3$-descent map is denoted 
$\wh{\al}$. Moreover we have $\wh{T}=(0,\wh{b}\sqrt{\wh{D}})=(0,(27b-4a^3D)\sqrt{-3D}/9)$.
\item The map $\phi$ from $E$ to $\wh{E}$ is defined by
$$\phi(P)=\left(\dfrac{x^3+4D((a^2/3)x^2+abx+b^2)}{x^2},\ \dfrac{y(x^3-4Db(ax+2b))}{x^3}\right)$$
for $P\ne O$ and $P\ne\pm T$, and $\phi(P)=\wh{O}$ if $P=O$ or $P=\pm T$,
and the map $\wh{\phi}$ from $\wh{E}$ to $E$ is defined in the same way,
replacing the coefficients of $E$ by those of $\wh{E}$, except that the 
$x$-coordinate must be divided by $9$ and the $y$-coordinate by $27$.
\end{enumerate}
\end{definition}

\begin{proposition}\begin{enumerate}\item $\phi$ and $\wh{\phi}$ are dual
$3$-isogenies (in particular group homomorphisms) between $E$ and $\wh{E}$,
so that $\wh{\phi}\circ\phi$ and $\phi\circ\wh{\phi}$ are the 
multiplication-by-$3$ maps on $E$ and $\wh{E}$ respectively. The kernel
of $\phi$ (over $\ov{\Q}$) is $\{O,\pm T\}$, and that of $\wh{\phi}$ is
$\{\wh{O},\pm\wh{T}\}$.
\item The map $\al$ is a group homomorphism from $E(\Q)$ to $G_3$,
and $\Ker(\al)=\Ima(\wh{\phi})$.
\end{enumerate}\end{proposition}

\section{$3$-Descent with a Rational $3$-Isogeny}

We now explain how the use of the $3$-descent maps $\al$ and $\wh{\al}$
gives a precise estimate on the rank of $E$ (and of the isogenous curve
$\wh{E}$, which has the same rank). Before proving the main result 
(Proposition \ref{alwh} below), we need the
following precise description of the rational $3$-torsion points of an
elliptic curve (evidently, if an elliptic curve does not have a rational
$3$-torsion subgroup, in other words if it does not have an equation of the
form $y^2=x^3+D(ax+b)^2$, the only rational $3$-torsion point is $O$).

\begin{lemma}\label{lemal} Let $y^2=x^3+D(ax+b)^2$ be the equation of an
elliptic curve $E$ with rational $3$-torsion subgroup, and assume as usual that
this equation is written so that $D$ is a fundamental discriminant. The
rational $3$-torsion points of $E$ are the following:
\begin{enumerate}\item If $D=1$, the points $O$ and $(0,\pm b)$.
\item If $D=-3$ and $2(9b+4a^3)=t^3$ is the cube of a rational number $t\neq 0$,
the point $O$ and the points $P$ such that $\displaystyle{x(P)=\frac{t^2}{3}+\frac{3}{t^2}\Big(4ab+\frac{16}{9}a^4\Big)+\frac{4a^2}{3}}$.
\item Otherwise, only the point $O$.
\end{enumerate}
\end{lemma}

\Proof Let $Q=(x,y)$ be a $3$-torsion point. Then we have $x([2]Q)=x(-Q)=x(Q)$, which gives, using the formulas on p. 59 of \cite{Sil1} for the duplication law on the elliptic curve $E$:
$$x(3x^{3}+4Da^{2}x^{2}+12Dabx+12Db^{2})=0.$$
Let $P(x)=3x^{3}+4Da^{2}x^{2}+12Dabx+12Db^{2}$. Note that $\mathrm{Disc}(P)=-48D^{2}(-27b+4Da^3)^2b^2$.
So either we have $x=0$, then $y^{2}=Db^2$, or we have $P(x)=0$ and we obtain after an easy calculation:
$$y^{2}=-\frac{D}{3}(ax+3b)^{2}.$$
It is then straightforward to find the rational solutions, keeping in mind that $D$ is a fundamental discriminant.
\fp

We can now give the following the exact analogue of Proposition 8.2.8 of
\cite{Cohen3}, whose proof we follow verbatim.

\begin{proposition}\label{alwh} Let $E$ be the elliptic curve
$y^2=x^3+D(ax+b)^2$ and $\wh{E}$ the $3$-isogenous curve with equation 
$y^2=x^3-3D(ax+(27b-4a^3D)/9)^2$ as above, and let $\al$ and $\wh{\al}$ be the
corresponding $3$-descent maps. Then
$$|\Ima(\al)||\Ima(\wh{\al})|=3^{r+\delta}\;,$$
where $r$ is the rank of $E$ (and of $\wh{E}$), and $\delta=1$ if $D=1$
or $D=-3$ and $\delta=0$ otherwise.\end{proposition}

\Proof If $E_t$ denotes the torsion subgroup of $E$ we have
$$E(\Q)/3E(\Q)\isom E_t(\Q)/3E_t(\Q)\oplus(\Z/3\Z)^r\;.$$
Set $G=E_t(\Q)$. We know that if $G$ is a finite Abelian group then
$G/3G$ is noncanonically isomorphic to $G[3]$, in other words to the group
of $3$-torsion points of $G$. Thus by Lemma \ref{lemal}, $E_t(\Q)/3E_t(\Q)$
is trivial unless either $D=1$, or $D=-3$ and $2(9b+4a^3)$ is a cube.
Write $\delta_{D,n}$ for the usual Kronecker $\delta$-symbol, and
$\ga(a,b)$ for the condition that $2(9b+4a^3)$ is a cube. With this notation
we can thus write
$$|E(\Q)/3E(\Q)|=3^{r+\delta_{D,1}+\delta_{D,-3}\ga(a,b)}\;.$$
On the other hand, let us consider our $3$-isogenies $\phi$ and $\wh{\phi}$.
Since $\wh{\phi}\circ\phi$ is the multiplication-by-$3$ map, we evidently have 
$$|E(\Q)/3E(\Q)|=[E(\Q):\wh{\phi}(\wh{E}(\Q))][\wh{\phi}(\wh{E}(\Q)):\wh{\phi}(\phi(E(\Q)))]\;.$$
Now for any group homomorphism $\wh{\phi}$ and subgroup $B$ of 
finite index in an abelian group $A$ we evidently have
$$\dfrac{\wh{\phi}(A)}{\wh{\phi}(B)}\isom \dfrac{A}{B+\Ker(\wh{\phi})}\isom \dfrac{A/B}{(B+\Ker(\wh{\phi}))/B}\isom\dfrac{A/B}{\Ker(\wh{\phi})/(\Ker(\wh{\phi})\cap B)}\;.$$
Thus 
$$[\wh{\phi}(A):\wh{\phi}(B)]=\dfrac{[A:B]}{[\Ker(\wh{\phi}):\Ker(\wh{\phi})\cap B]}\;.$$
We are going to use this formula with $A=\wh{E}(\Q)$ and $B=\phi(E(\Q))$.
We know that $\Ker(\wh{\phi})$ (over $\ov{\Q}$) has three elements $\wh{O}$ 
and $\pm\wh{T}$, and $\wh{T}\in\phi(E(\Q))$ if and only if $D=-3$, so (once
again over $\ov{\Q}$), 
$\Ker(\wh{\phi})=\{O,\pm\wh{T}\}$ if $D=-3$, and is trivial otherwise.
Thus, if $D\ne-3$ we have $[\Ker(\wh{\phi}):\Ker(\wh{\phi})\cap B]=1$.
Assume now that $D=-3$, so that the equation of $E$ is $y^2=x^3-3(ax+b)^2$,
and that of $\wh{E}$ can be taken to be $y^2=x^3+(3ax+(9b+4a^3))^2$.
Then $[\Ker(\wh{\phi}):\Ker(\wh{\phi})\cap B]=1$ if
$\wh{T}\in\phi(E(\Q))$, and is equal to $3$ otherwise. We know (see
for instance \cite{Cohen3}, Proposition 8.4.4) that $\wh{T}\in\phi(E(\Q))$
if and only if $2(9b+4a^3)$ is a cube, in other words with the notation
introduced above, if and only if $\ga(a,b)=1$. Thus,
$$[\wh{\phi}(\wh{E}(\Q)):\wh{\phi}(\phi(E(\Q)))]=\dfrac{[\wh{E}(\Q):\phi(E(\Q))]}{3^{\delta_{D,-3}(1-\ga(a,b))}}\;.$$
Putting everything together we obtain
\begin{align*}3^{r+\delta_{D,1}+\delta_{D,-3}\ga(a,b)}&=|E(\Q)/3E(\Q)|=[E(\Q):\wh{\phi}(\wh{E}(\Q))]
[\wh{\phi}(\wh{E}(\Q)):\wh{\phi}(\phi(E(\Q)))]\\
&=[E(\Q):\wh{\phi}(\wh{E}(\Q))][\wh{E}(\Q):\phi(E(\Q))]3^{\delta_{D,-3}(\ga(a,b)-1)}\;.\end{align*}
On the other hand the $3$-descent map $\al$ on $E(\Q)$ has kernel 
$\wh{\phi}(\wh{E}(\Q))$, so $[E(\Q):\wh{\phi}(\wh{E}(\Q))]=|\Ima(\al)|$, and
similarly $[\wh{E}(\Q):\phi(E(\Q))]=|\Ima(\wh{\al})|$, so finally we obtain
$$|\Ima(\al)||\Ima(\wh{\al})|=3^{r+\delta_{D,1}+\delta_{D,-3}}\;,$$
proving the proposition.\fp

It follows from this proposition that to compute the rank it is sufficient
to compute the cardinality of $\Ima(\al)$ and of $\Ima(\wh{\al})$, which
we do separately. As in the case of $2$-descent, we cannot give an
algorithm for this, since there is an obstruction embodied in the $3$-part
of the Tate--Shafarevich group of $E$, but the method works in many cases.
The goal of the next sections is thus to compute $|\Ima(\al)|$.

\section{The Case $D=1$}

We first treat the case $D=1$.
We choose the equation of our elliptic curve as $y^2=x^3+(ax+b)^2$ with
$a$ and $b$ as in Lemma \ref{lemeq}, and we recall that the fundamental
$3$-descent map $\al$ from $E(\Q)$ to $\Q^*/{\Q^*}^3$ is defined by 
$\al(O)=1$, $\al((0,b))=1/(2b)$, and $\al((x,y))=y-(ax+b)$ for all other 
points of $E(\Q)$.

\smallskip

{\bf Note:} To avoid unnecessary wording, when we speak of a solution to
a homogeneous equation we always mean a \emph{nontrivial} solution, where the
variables are not all equal to $0$. Similarly, when we speak of a solution to
a homogeneous congruence modulo $p^k$ for some prime $p$, we always mean a
solution where all the variables are $p$-integral, and at least one of them
having $p$-adic valuation equal to $0$, typically 
$\min(v_p(X),v_p(Y),v_p(Z))=0$.

\begin{theorem}\label{thd1} Keep the above notation. 
\begin{enumerate}\item An element $\ov{u}\in\Q^*/{\Q^*}^3$
belongs to the image of $\al$ if and only if for some (or any) representative
$u\in\Q^*$ the homogeneous cubic equation
$$uX^3+(1/u)Y^3+2bZ^3-2aXYZ=0$$
has an integer (or rational) solution.
\item More precisely, for $u=1$ it has the solution $(X,Y,Z)=(1,-1,0)$,
for $u=1/(2b)$ it has the solution $(X,Y,Z)=(0,1,-1)$, and if 
$y-(ax+b)=uz^3$ for some $z\in\Q^*$ it has the solution 
$(X,Y,Z)=(z^2,-x,z)$. Conversely, if $(X,Y,Z)$ is a solution of the cubic
with $Z\ne0$ then $(x,y)=(-XY/Z^2,(uX^3-(1/u)Y^3)/(2Z^3))$ is a preimage of 
$u$ in $E(\Q)$, and $z=X/Z$.
\item If the above cubic has a rational solution and if
$u$ is the unique positive integer cubefree representative of $\ov{u}$, then
$u_1u_2\mid(2b)$, where $u=u_1^2u_2$ with the $u_i$ squarefree and coprime, and
the solubility of the cubic is equivalent to that of
$$u_1X^3+u_2Y^3+(2b/(u_1u_2))Z^3-2aXYZ=0\;.$$
\end{enumerate}
\end{theorem}
\Proof (1) and (2). The cases $u=1$ and $u=1/(2b)$ (corresponding to the 
points $O$ and $T=(0,b)$ respectively) being clear, we assume that we are not
in these cases. Then by definition if $\ov{u}$ belongs to the image of $\al$ 
there exists $(x,y)\in E(\Q)$ and $z\in\Q^*$ such that $uz^3=y-(ax+b)$, and 
if we set $X=z^2$, $Y=-x$, and $Z=z$ then
\begin{align*}uX^3+(1/u)Y^3+2bZ^3-2aXYZ&=\dfrac{1}{u}(u^2z^6+2uz^3(ax+b)-x^3)\\
&=\dfrac{1}{u}((uz^3+ax+b)^2-x^3-(ax+b)^2)\\
&=\dfrac{1}{u}(y^2-(x^3+(ax+b)^2))=0\;,\end{align*}
as claimed. Note that since $z\in\Q^*$ we have $Z\ne0$. Conversely, let
$(X,Y,Z)$ be a solution to our cubic with $Z\ne0$. If we set
$x=-XY/Z^2$ and $y=(uX^3-(1/u)Y^3)/(2Z^3)$, we have
$x^3+(ax+b)^2=(-X^3Y^3+Z^2(bZ^2-aXY)^2)/Z^6$, and since by the cubic equation
we have $-2Z(bZ^2-aXY)=uX^3+(1/u)Y^3$, it follows that
\begin{align*}x^3+(ax+b)^2&=((uX^3+(1/u)Y^3)^2-4X^3Y^3)/(4Z^6)\\
&=(uX^3-(1/u)Y^3)/(4Z^6)=y^2\;,\end{align*}
so $(x,y)\in E(\Q)$. Furthermore, 
\begin{align*}\al((x,y))&=y-(ax+b)=(uX^3-(1/u)Y^3)/(2Z^3)-(bZ^2-aXY)/Z^2\\
&=(1/(2Z^3))(uX^3-(1/u)Y^3-2Z(bZ^2-aXY))\\
&=(1/(2Z^3))(2uX^3)=u(X/Z)^3\;,\end{align*}
so is equal to $u$ up to cubes, hence $(x,y)$ is indeed a preimage of $u$,
as claimed.

\smallskip

For (3), let $u$ be the (unique) positive integer cubefree representative of $\ov{u}$,
and write $u=u_1^2u_2$ with the $u_i$ squarefree and coprime. Replacing
$Y$ by $u_1u_2Y$ in the cubic equation we obtain the equivalent equation
$$u_1^2u_2X^3+u_1u_2^2Y^3+2bZ^3-2au_1u_2XYZ=0\;.$$
It is clear that this homogenous cubic has a rational solution if 
and only if it has an integer solution, and we may in addition 
assume that $\gcd(X,Y,Z)=1$. Assume by contradiction that $u_1u_2\nmid 2b$.
Since $u_1u_2$ is squarefree, this means that there exists a prime $p$ such
that $p\mid u_1u_2$ and $p\nmid 2b$. Since exchanging $X$ and $Y$ in the
above equation is equivalent to the exchange of $u_1$ and $u_2$, we may
assume that $p\mid u_1$, hence $p\nmid u_2$. Since $p$ divides the first,
second and fourth term of the equation it divides the third, and since 
$p\nmid 2b$, we deduce that $p\mid Z$. Thus $p^2$ divides the first,
third and fourth term, so it divides the second, and since $u_1$ is squarefree
and $p\nmid u_2$ we deduce that $p\mid Y$. Thus, $p^3$ divides the
second, third, and fourth term, so it divides the first, and again since
$u_1$ is squarefree and $p\nmid u_2$, we deduce that $p\mid X$, 
contradicting the assumption $\gcd(X,Y,Z)=1$. We can thus divide by $u_1u_2$
to obtain the final equation given in (3).\fp

\textit{Geometric Interpretation:}

The plane cubic given in (1) of Theorem \ref{thd1} is the equation of a twist $C$ of the elliptic 
curve $E$. We can recover a linear action by linear automorphisms by doing
the following: we pick $s\in{\mathbb{Q}^{*}}$ and consider the action: 
$$(X:Y:Z)\LR (sX:(1/s)Y:Z).$$
This action gives a curve isomorphic to $C$ with $u$ replaced by $us^{3}$. 
This is what was predicted in the geometrical interpretation of 
\cite{CreFiOSiSto1,CreFiOSiSto2} recalled in the introduction.

\vspace{0.5cm}

\section{The Case $D\ne1$}

We now assume specifically that $D\ne1$, so that $K=\Q(\sqrt{D})$ is a
genuine quadratic field. Note that to use $3$-descent with a rational
$3$-torsion subgroup, we must compute the image of the $3$-descent map
both for the curve $E$, and for a $3$-isogenous curve $\wh{E}$, whose
$\wh{D}$ is such that $\wh{D}=-3D$. Thus, we always need to
consider curves with $D\ne1$. We will denote by $\tau$ the nontrivial element 
of the Galois group $\Gal(K/\Q)$, so that $\tau(\sqrt{D})=-\sqrt{D}$,
and by $\mathrm{N}=\mathrm{N}_{K/\Q}$ the norm from $K$ to $\Q$, both for elements and for
ideals. If $u\in K^*$ we denote by $[u]$ the class of $u$ in $K^*/{K^*}^3$.

\subsection{Description of the Image of $\al$}

The equation of our curve is $y^2=x^3+D(ax+b)^2$, and the $3$-descent map is a
map $\al$ from $E(\Q)$ to the subgroup $G_3$ of $K^*/{K^*}^3$ of classes $[u]$
of elements $u$ such that $u\tau(u)=\mathrm{N}_{K/\Q}(u)\in{\Q^*}^3$, defined by
$\al(O)=1$ and $\al((x,y))=y-(ax+b)\sqrt{D}$ for all other points of
$E(\Q)$ (note that $T=(0,b\sqrt{D})\notin E(\Q)$). The image of $\al$
can be described as follows.

\begin{theorem}\label{thdnot1} Keep the above notation. 
\begin{enumerate}\item An element $[u]\in G_3\subset K^*/{K^*}^3$
belongs to the image of $\al$ if and only if for some (or any) representative
$u\in K^*$ of the form $u=v^2\tau(v)$ the homogeneous cubic equation
$$2v_2X^3+2Dv_1Y^3+\dfrac{2b}{v_1^2-Dv_2^2}Z^3+6v_1X^2Y+6v_2DXY^2+2a(X^2Z-DY^2Z)=0$$
has an integer (or a rational) solution, where we write
$v=v_1+v_2\sqrt{D}$.
\item More precisely, for $u=1$ it has the solution $(X,Y,Z)=(1,0,0)$, and if 
$y-(ax+b)\sqrt{D}=v^2\tau(v)z^3$ for some $z\in\Q^*$ it has the solution 
$(X,Y,Z)=(z_1,z_2,1)$, where $z=z_1+z_2\sqrt{D}$.
Conversely, if $(X,Y,Z)$ is a solution of the cubic with $Z\ne0$ then 
$$(x,y)=\left(v\tau(v)\dfrac{X^2-DY^2}{Z^2},v\tau(v)\dfrac{\Re(v(X+Y\sqrt{D})^3)}{Z^3}\right)$$
is a preimage of $u$ in $E(\Q)$ with $z=(X+Y\sqrt{D})/Z$, where by abuse of 
notation we write $\Re(\al)=(\al+\tau(\al))/2$.
\end{enumerate}
\end{theorem}

\Proof If we set $a'=a\sqrt{D}$ and $b'=b\sqrt{D}$, and ignore for the moment
all rationality questions, the equation of our curve is $y^2=x^3+(a'x+b')^2$,
so if $y-(a'x+b')=y-(ax+b)\sqrt{D}=uz^3$ then $(X,Y,Z)=(z^2,-x,z)$ is a 
solution to the modified cubic 
$uX^3+(1/u)Y^3+2b\sqrt{D}Z^3-2a\sqrt{D}XYZ=0$, and conversely if
$(X,Y,Z)$ is such a solution then $(x,y)=(-XY/Z^2,(uX^3-(1/u)Y^3)/(2Z^3))$ is
a point on the curve. So formally there is no problem. We now must
add the condition that not only $(x,y)\in E(K)$, but that $(x,y)\in E(\Q)$.

We first choose suitable representatives $u$. Let for the moment $u$ be
any representative of $[u]$, so that $u\tau(u)=\mu^3$ for some $\mu\in\Q$. 
This implies that $(u^2/\mu^3)\tau(u^2/\mu^3)=1$, hence by a weak form of
Hilbert's theorem 90, we have $u^2=\mu^3 v/\tau(v)$ for some $v\in K$, hence 
$u=(u/\mu)^3\tau(v)/v=(u/(\mu v))^3 v^2\tau(v)$. Since $u$ is defined up to 
cubes, we may therefore assume that $u=v^2\tau(v)$ for some $v\in K^*$. Thus 
$u\tau(u)=\mu^3$ with $\mu=v\tau(v)$. Multiplying $v$ by a suitable rational
number we can assume that $v\in\Z_K$.

Now the condition $x\in \Q$ means that $XY/Z^2\in\Q$, so that 
$Y/Z=\la\tau(X/Z)$ for some $\la\in\Q^*$. The condition $y\in \Q$ is
thus that $u\al^3-(\la^3/u)\tau(\al)^3\in\Q$, where $\al=X/Z$.
Replacing $u$ by $v^2\tau(v)$ and dividing by the rational number $v\tau(v)$
gives the condition $v\al^3-(\la^3/(v^3\tau(v)^2))\tau(\al)^3\in\Q$.
Setting $\be=v\al^3$ and $r=(\la/(v\tau(v)))^3\in\Q^*$, this gives the
condition $\be-r\tau(\be)\in\Q$, so if $\be=s+t\sqrt{D}$, we have 
$\be-r\tau(\be)=s+t\sqrt{D}-r(s-t\sqrt{D})$, hence the condition is $t(r+1)=0$,
in other words either $r=-1$, so $\la=-\mu$, or $t=0$. 

If $t=0$ then $\be=v\al^3\in\Q$, hence $u=v^2\tau(v)=(\be/(\al^2\tau(\al)))^3$,
so $[u]$ is trivial. It follows that the case $t=0$ corresponds to
the unit element of $G_3$, which we now exclude.

Thus we may assume that $\la=-\mu$. The cubic equation is thus
$u\al^3-\tau(u\al^3)+2a\sqrt{D}\mu\al\tau(\al)+2b\sqrt{D}=0$,
and since $\mu=v\tau(v)$ this gives
$(v\al^3-\tau(v\al^3))/(2\sqrt{D})+a\al\tau(\al)+b/(v\tau(v))=0$.
Recall that $v$ is given, so write $v=v_1+v_2\sqrt{D}$ and
$\al=x_1+x_2\sqrt{D}$. The above equation is thus
$$2v_2x_1^3+2Dv_1x_2^3+6v_1x_1^2x_2+6v_2Dx_1x_2^2+2a(x_1^2-Dx_2^2)+2b/(v_1^2-Dv_2^2)=0\;,$$
so setting $x_1=X/Z$, $x_2=Y/Z$ with $X$, $Y$, $Z$ in $\Z$ we obtain finally
$$2v_2X^3+2Dv_1Y^3+(2b/(v_1^2-Dv_2^2))Z^3+6v_1X^2Y+6v_2DXY^2+2a(X^2Z-DY^2Z)=0\;.$$
The formulas for $(X,Y,Z)$ knowing $x$ and $y$ and vice versa are immediately
obtained by replacing the corresponding quantities in the formula given for
the case $D=1$.\fp

\smallskip

\noindent
{\bf Remarks.}\begin{enumerate}\item The cubic equation in (1) of Theorem \ref{thdnot1} can be written
\begin{align*}&(v(X+Y\sqrt{D})^3-\tau(v)(X-Y\sqrt{D})^3)/\sqrt{D}\\
&\phantom{=}+2aZ(X+Y\sqrt{D})(X-Y\sqrt{D})+2b/(v\tau(v))Z^3=0\;.\end{align*}
\item By the theorem, the solubility of the equation depends only on the class
$[u]$ of $u$, hence we can change $v$ into $v\ga^3$ for any
$\ga\in K^*$, or $v$ into $vr$ for any $r\in\Q^*$ without changing the
solubility of the equation. This is of course clear directly.
\item Since the image $\Ima(\al)$ of $\al$ is a \emph{group}, 
$[u]\in\Ima(\al)$ 
if and only if $[1/u]\in\Ima(\al)$, so the solubility for $v$ is equivalent to 
that for $v^{-1}$. Furthermore since $\tau(u)=u^{-1}(v\tau(v))^3$, 
$[\tau(u)]\in \Ima(\al)$, so the solubility for $v$ is equivalent to that for 
$\tau(v)$.
\item \textit{Geometric Interpretation:} The plane cubic given in (1) of Theorem \ref{thdnot1} is the 
equation of a twist $C$ of the elliptic curve $E$. We get a linear action by
linear automorphisms by doing the following: we pick $s\in{\mathbb{Q}^{*}}$ 
and consider the action: $$(X:Y:Z)\LR (sX:sY:(1/s^{2})Z).$$
This action gives a curve isomorphic to $C$ with $u$ replaced by $us^{9}$.
\end{enumerate}

\smallskip

The reader will notice that we have not given an analogous result to (3) of
Theorem \ref{thd1}, which is essential since it is necessary to check only
a finite number of elements of $G_3$. We do this in the next subsection.

\subsection{Reduction of Elements of $G_3$}\label{secg3}

We begin by the following lemma.

\begin{lemma}\label{lemgcd} Assume that $x$ and $y$ are rational numbers such
that $y^2=x^3+D(ax+b)^2$, and write $x=m/d^2$ and $y=n/d^3$ with
$\gcd(m,d)=\gcd(n,d)=1$ and $d>0$. Finally, set 
$\f=\gcd(n-d(am+bd^2)\sqrt{D},n+d(am+bd^2)\sqrt{D})$, where the GCD is 
understood in the sense of ideals.
\begin{enumerate}\item There exist integers $f$ and $g$ such that $g$ is
a squarefree integer dividing $D$, and $\f=fg\gd$, where $\gd$ is the unique
ideal such that $\gd^2=g\Z_K$ (when $D=1$ we have of course $g=1$ and 
$\gd=\Z$).
\item If we write $f=f_1q^3$ with $f_1$ cubefree then $gf_1q^2\mid 2b$, and
in particular, $g\mid\gcd(D,2b)$.
\item If $p\mid f_1$ then $p$ is split in $K/\Q$, and in particular $g$ and
$f_1$ are coprime.
\end{enumerate}\end{lemma}

\Proof The case $D=1$, which we do not need, is left to the reader, so
assume that $D\ne1$, so that $K$ is a quadratic field. We can write uniquely 
$\f=F\gd$, where $F\in\Z$ and $\gd$ is an integral ideal of $K$ which is 
\emph{primitive}, in other words which is not divisible by any element of $\Z$
other than $\pm1$. Evidently $\gd$ cannot be divisible by inert primes; since
$\f$ is the GCD of two conjugate elements it is stable by conjugation,
hence $\gd$ cannot be divisible by an ideal $\p$ above a split prime $p$,
otherwise it would also be divisible by $\tau(\p)$, hence by $p=\p\tau(\p)$.
Finally since $\p^2=p\Z_K$ for a ramified prime $p$, $\gd$ cannot be divisible
by a ramified prime to a power higher than the first, so $\gd$ is equal to
a product of distinct ramified primes. Thus we have $\gd=\prod_{\p\in S_{0}}\p$ for
some set $S_{0}$ of ramified primes $\p$, and if $g=\prod_{\p\in S_{0}}p$, where $p$
is the prime number below $\p$, we thus have $\f=F\prod_{p\mid g}\p$ and
$g\mid D$, hence also $\f^2=F^2g$.

\smallskip

Let us now use our equation. Replacing $x$ and $y$ by $m/d^2$ and $n/d^3$ we
obtain the equation 
$$(n-d(am+bd^2)\sqrt{D})(n+d(am+bd^2)\sqrt{D})=n^2-Dd^2(am+bd^2)^2=m^3\;.$$
It follows that $F^2g\mid m^3$. Let $p\mid g$ and $\p$ be the prime ideal
above $p$. Since the two factors are conjugate, if $v\ge1$ is the $\p$-adic
valuation of the first factor, it is also that of the second. This implies
both that $v=v_{\p}(\f)$ and that $3v_{\p}(m)=2v_{\p}(\f)$, in other words
$$3v_p(m)=v_{\p}(\f)=v=v_{\p}(\f^2)/2=v_{\p}(F^2g)/2=v_p(F^2g)=1+2v_p(F)$$
since $p\mid g$ and $g$ is squarefree. We deduce that $v_p(F)\equiv1\pmod3$,
and in particular that $v_p(F)\ge1$, so $p\mid F$, proving that $g\mid F$.
Thus, we will set $F=fg$ for some $f\in\Z$. The same reasoning shows that if 
$p$ is \emph{any} inert or ramified prime (dividing $\f$ or not) then 
$3v_p(m)=v_{\p}(\f)$, so $3\mid v_{\p}(\f)$.

Since $f^2g^3\mid m^3$ we have $g\mid m$ and $f^2\mid (m/g)^3$.
Write $f=f_1q^3$ with $f_1$ cubefree. For all primes $p$ we have 
$v_p(m/g)\ge v_p(q)+\lceil 2v_p(f_1)/3\rceil$. Since $0\le v_p(f_1)\le2$
we have $\lceil 2v_p(f_1)/3\rceil=v_p(f_1)$, so $f_1q^2\mid m/g$.
Note that $\f^2=f^2g^3=f_1^2g^3q^6$, so since
for any inert or ramified prime we have $3\mid v_{\p}(\f)$, for such a prime
we have $v_p(f_1)=0$, so $f_1$ is only divisible by split primes. In particular
it is coprime to $D$, hence to $g$.

Since $(n-d(am+bd^2)\sqrt{D})/(fg)$ is an algebraic integer and $D$ is a 
fundamental discriminant, it follows that $fg\mid2\gcd(n,d(am+bd^2))$, hence
$fg=gf_1q^3\mid 2\gcd(n,d(am+bd^2))$. Since $2bd^3=2d(am+bd^2)-2adm$ and
$gf_1q^2\mid m$ we deduce that $2bd^3\equiv0\pmod{gf_1q^2}$. Since $d$ and
$n$ are coprime there exist integers $u$ and $v$ such that $un+vd^3=1$, hence
$2b=2bvd^3+2bun$, and since $gf_1q^2\mid 2bd^3$ and $gf_1q^2\mid 2n$, we have
$gf_1q^2\mid 2b$, as claimed.\fp

\begin{corollary}\label{corq3} Keep the above notation, and let
$[u]\in\Ima(\al)\subset G_3$. There exists an integral ideal $\v$ of $K$
such that $u\Z_K=\v^2\tau(\v)\q^3$ for some ideal $\q$ of $K$, and which is 
such that $\gcd(\v,\tau(\v))=1$ and $f_1=\mathrm{N}_{K/\Q}(\v)$ is a cubefree divisor
of $2b$ divisible only by primes which are split in $K/\Q$.
\end{corollary}

\Proof The above lemma states that if we set
$$\f=\gcd\Big(n-d(am+bd^2)\sqrt{D},n+d(am+bd^2)\sqrt{D}\Big)$$
there exist integers $f_1$, $g$, and $q$ in $\Z$ and an ideal $\gd\in K$ such
that $\f=f_1q^3g\gd$ with $f_1$ cubefree divisible only by split primes, 
$g\mid \gcd(D,2b)$, $g\Z_K=\gd^2$, and $gf_1q^2\mid 2b$.
Thus $\f=f_1(q\gd)^3$. Set $\a_-=(n-d(am+bd^2)\sqrt{D})/\f$
and $\a_+=(n+d(am+bd^2)\sqrt{D})/\f$, so that $\gcd(\a_-,\a_+)=1$
and our equation implies $\a_-\a_+=m^3/\f^2$ (since we work with ideals,
we lose some unit information here). Since $f_1$ is cubefree and divisible
only by split primes, it is also cubefree as an ideal, and the condition 
$\f^2\mid m^3$ implies as above that $f_1(q\gd)^2\mid m$, hence 
$gf_1q^2\mid m$ (this time in $\Z$), so the equation now reads
$$\a_-\a_+=f_1(m/(gf_1q^2))^3\Z_K\;.$$
Since $\f$ is stable by conjugation, it is clear that $\a_+=\tau(\a_-)$ so
this equation gives the norm of $\a_{\pm}$. In any case, write 
$\a_-=\f_-\q_-^3$, $\a_+=\f_+\q_+^3$ with $\f_{\pm}$ cubefree, so that in 
particular $\gcd(f_-,f_+)=1$, and also $\f_+=\tau(\f_-)$ and
$\q_+=\tau(\q_-)$. At the level of ideals our equation is thus
$\f_-\f_+(\q_-\q_+)^3=f_1(m/(gf_1q^2))^3\Z_K$. As we already mentioned
$f_1$ is also cubefree in $K$, and since $\f_-$ and $\f_+$ are coprime and 
cubefree, by uniqueness of the decomposition as the product of a cubefree 
ideal with a cube it follows that $\q_-\q_+=m/(gf_1q^2)\Z_K$ and
$\f_-\f_+=f_1\Z_K$. In particular, $\f_-\mid f_1$.

Now recall that the $3$-descent map $\al$ is defined on an affine point
as the class modulo cubes of $y-(ax+b)\sqrt{D}=(n-d(am+bd^2)\sqrt{D})/d^3$.
Thus
$$\al((x,y))\Z_K=\f\a_-=f_1(q\gd)^3\f_-\q_-^3=f_1\f_-(q\gd\q_-)^3=\v^2\tau(\v)\q_1^3$$
for some ideal $\q_1$, with $\v=\f_-$, as claimed.\fp

\smallskip

\noindent
{\bf Remark.} It is clear that the condition $\gcd(\v,\tau(\v))=1$ implies that
$\v$ is primitive, in other words that the only elements of $\Z$ which divide
it are $\pm1$. Furthermore, since $\v^2\tau(\v)=u\q^{-3}$, the ideal class of
$\v$ in $Cl(K)$ belongs in fact to $(Cl(K)/Cl(K)^3)[\tau+2]$, where for any
group $G$ and map $\phi$ from $G$ to $G$, $G[\phi]$ denotes the elements of
$G$ killed by $\phi$, in other words the kernel of $\phi$.

\smallskip

Now recall the definition of a \emph{$3$-virtual unit} (or \emph{virtual cube}) and \emph{$3$-Selmer group}: an 
element $u\in K^*$ is a virtual cube if $u\Z_K=\q^3$ is the cube of an ideal. 
The group of virtual cubes modulo cubes of elements is called the $3$-Selmer
group of $K$ and denoted $S_3(K)$. It is clear that $S_3(K)\subset G_3$, and 
it is well-known and easy that we have a natural exact sequence
$$1\LR U(K)/U(K)^3\LR S_3(K)\LR Cl(K)[3]\LR 1\;.$$
Let as usual $I(K)$ denote the group of (nonzero) fractional ideals of $K$,
and let $G_3^i$ be the subgroup of $I(K)/I(K)^3$ of classes of ideals whose 
norm is a cube.

\begin{lemma} We have a natural exact sequence
$$1\LR S_3(K)\LR G_3\LR G_3^i\LR Cl(K)/Cl(K)^3\LR1\;.$$
\end{lemma}

\Proof Consider first the map $i$ from $G_3$ to $G_3^i$ which sends a class
$[u]$ to the class of $u\Z_K$. It is clear that it does send $G_3$ to $G_3^i$.
If $[u]$ is sent to the unit element of
$G_3^i$ this means that $u\Z_K=\q^3$ for some ideal $\q$, in other words
that $[u]\in S_3(K)$, giving the kernel. Consider now the
map sending the class of an ideal modulo cubes to its ideal class. It defines
a map $\pi$ from $G_3^i$ to $Cl(K)/Cl(K)^3$. If some ideal $\a$ is sent
to the unit element of $Cl(K)/Cl(K)^3$ this means that there exist an ideal $\q$
and an element $\ga\in K^*$ such that $\a=\ga\q^3$. Thus, the class of $\a$
modulo cubes of ideals is equal to that of $\ga\Z_K$. Furthermore, the norm
of $\a$ is a cube, so that of $\ga$ also (since $-1$ is a cube), hence $\ga$
does come from $G_3$, proving exactness at $G_3^i$. Finally, let us show that
the map $\pi$ is surjective. Let $\a$ be an ideal, representative of an 
element of $Cl(K)/Cl(K)^3$. Then since $\mathrm{N}_{K/\Q}(\a)=a\in\Q^*$,
$\a\mathrm{N}_{K/\Q}(\a)$ is in the same ideal class as $\a$, and its norm is evidently
equal to $(\mathrm{N}_{K/\Q}(\a))^3$, so the class of $\a$ is the image of the
class of $\a\mathrm{N}_{K/\Q}(\a)$ in $G_3^i$, proving surjectivity and the lemma.\fp

\begin{corollary}\label{corali} If we denote by $\al^i=i\circ\al$ the 
$3$-descent map from $E$ to $G_3^i$ we have
$$|\Ima(\al)|=|\Ima(\al^i)||S_3(K)\cap\Ima(\al)|\;.$$
In particular, if $D<0$, $D\ne-3$, and $3\nmid h(K)$ we have
$|\Ima(\al)|=|\Ima(\al^i)|$.
\end{corollary}

\Proof By the above lemma the map $i$ induces an injection from 
$G_3/S_3(K)$ to $G_3^i$. Thus, for any subgroup $H$ of $G_3$ the map $i$
induces a bijection from $H/(S_3(K)\cap H)$ to $i(H)$. Applying this to the
finite subgroup $\Ima(\al)$ gives the formula of the corollary, and the
special case corresponds to $S_3(K)=1$.\fp

\smallskip

Note that if $\ga$ is a $3$-virtual unit then $\mathrm{N}(\ga)=\ga\tau(\ga)$ is a cube.
More generally, if $\ga$ is such that $\mathrm{N}(\ga)=n^3$ is a cube we can write
$\ga=v^2\tau(v)q^3$ with $v=\ga$ and $q=1/n$. Thus, using Corollaries
\ref{corq3} and \ref{corali} and Theorem \ref{thdnot1} we can compute
$|\Ima(\al)|$. Note however that it will be completely algorithmic and easy
to prove everywhere local solubility of our homogeneous cubics, so that
we will compute the $3$-Selmer group of $E$. On the other hand, for global
solubility, either we find a solution with a reasonable search bound, or
we are led to believe that no such solution exists, coming from a nontrivial
element of the $3$-part of the Tate--Shafarevich group. Thus, it is reasonable
to give an algorithm which computes both the rank of the $3$-Selmer group
and a lower bound for the rank of the curve, so which gives the exact rank
when they coincide.

\begin{enumerate}\item For each class $[\ga]\in S_3(K)$ choose a representative
$\ga\in K^*$, and check that the cubic of Theorem \ref{thdnot1} has a solution
for $v=\ga$ (more on this later), let $T_S$ be the group of $[\ga]$ for which 
it is everywhere locally soluble (ELS), and $T_G$ the subgroup generated by the elements for which we find 
that it has a global solution. Thus $T_G\subset T\subset T_S\subset S_3(K)$,
where $T=S_3(K)\cap\Ima(\al)$. This will allow us to compute $T$ exactly
if $T_S=T_G$, and otherwise if the search bound is sufficiently large,
we suspect (but cannot prove without further work) that in fact $T=T_G$ and
that the elements of $T_S/T_G$ correspond to nontrivial elements of the
$3$-part of the Tate--Shafarevich group. Note that using the fact that $T_G$ 
and $T_S$ are groups, it is not necessary to test all classes $[\ga]$, but in 
fact using the fact that they are even $\F_3$-vector spaces it is sufficient 
to work on bases of these spaces and use linear algebra. Finally, choose a set
$R_S$ of representatives of $S_3(K)/T_S$ and a set $R_G$ of representatives of 
$T_S/T_G$.
\item Let $f$ be the largest positive integer cubefree divisor of $2b$ divisible only
by split primes, write $f=\prod_{1\le i\le s}p_i^{v_i}$ with $1\le v_i\le 2$.
For each $p_i$, let $\p_i$ be one of the two prime ideals above $p_i$,
fixed once and for all. Find all ideals $\v$ of the form
$\v=\prod_{1\le i\le s}\p_i^{x_iv_i}$ with $-1\le x_i\le 1$ whose ideal
class is a cube (although this seems to be $3^s$ principal ideal tests,
it is easy to reduce to only $s$ such tests using linear algebra, but in
practice $Cl(K)$ will be small).
\item For each ideal $\v$ that we have found write $\v=u\q^3$ for some ideal 
$\q$ and some element $u\in K^*$, where clearly the class $[u]$ of $u$ in 
$G_3$ is determined uniquely modulo multiplication by an element of $S_3(K)$. 
For the moment, we choose any $u$ as above.
\item For each $[\ga_S]\in R_S$, where $R_S$ is the system of representatives
of $S_3(K)/T_S$ computed in (1), and any representative $\ga_S\in K^*$, check 
whether the cubic of Theorem \ref{thdnot1} is everywhere locally soluble (ELS) for $v=u\ga_S$ (more on this
later). If for some $[\ga_S]\in R_S$ this is the case, there will be only
one, and then $\v\in \Sel(\al^i)$, with evident notation, otherwise
$\v\notin \Sel(\al^i)$. In this latter case, we do nothing more, otherwise
for each $[\ga_G]\in R_G$, where $R_G$ is the system of representatives
of $T_S/T_G$ computed in (1), and any representative $\ga_G\in K^*$, check
whether the cubic of Theorem \ref{thdnot1} has a global solution up to a
reasonable search bound for $v=u\ga_S\ga_R$. If this is the case, once again
there will be only one, and then $\v\in\Ima(\al^i)$, otherwise we suspect
(but cannot be sure) that $\v\notin\Ima(\al^i)$. We let $I_G$ be the group
generated by the $\v$ for which we are sure, so that $I_G$ is a subgroup of
$\Ima(\al^i)$, probably equal to it.
\item At the end of this process we have computed the Selmer group cardinality
$|\Sel(\al)|=|\Sel(\al^i)||T_S|$, and groups $|I_G|$ and $T_G$ probably
equal to $\Ima(\al^i)$ and $T$ respectively, but in any case satisfying
$T_G\subset T\subset T_S$ and $I_G\subset\Ima(\al^i)\subset\Sel(\al^i)$,
and so $|I_G||T_G|\mid |\Ima(\al)|\mid |\Sel(\al)|$, the unknown quantity
being $|\Ima(\al)|$. If $|I_G||T_G|=|\Sel(\al)|$ then of course these
quantities are equal to $|\Ima(\al)|$. Otherwise, we output both quantities,
and a message saying that we expect $|\Ima(\al)|$ to be equal to
$|I_G||T_G|$ and that the elements of $\Sel(\al)/I_GT_G$ correspond to 
nontrivial elements of the $3$-part of the Tate--Shafarevich group.
\end{enumerate}

\section{Local Solubility of $u_1X^3+u_2Y^3+u_3Z^3-cXYZ=0$}

There remains to decide whether or not the cubics of Theorems \ref{thd1} and
\ref{thdnot1} have a rational solution. It is unfortunately
well-known that there is no algorithm for doing this. We thus proceed
as follows: we first check whether the cubic is everywhere locally soluble
(which we will abbreviate to ELS). If not, there are no rational solutions. 
Otherwise, either there is an obstruction in the $3$-part of the
Tate--Shafarevich group, or there does exist a rational
solution which we can find using a more or less efficient search. If we do
find one, we are done, otherwise we give up and can only give bounds on
$|\Ima(\al)|$, not its precise value.

Testing ELS is an algorithmic process, but is not
completely trivial. In this section we give such an algorithm. Since the 
degree is odd there is no need to look at local solubility in $\R$. 
We treat the following slightly more general problem: decide solubility in
$\Q_p$ of the equation
$$F(X,Y,Z)=u_1X^3+u_2Y^3+u_3Z^3-cXYZ=0\;.$$
For $D=1$, in other words, for the cubic of Theorem \ref{thd1},
we apply the results that we obtain with $u_3=2b/(u_1u_2)$ and $c=2a$.

\subsection{Reduction of the Cubic and Bad Primes}

By multiplying $X$, $Y$, and $Z$ by suitable powers of $p$ it is clear that
without loss of generality we may assume that $u_1$, $u_2$, $u_3$ and $c$
are $p$-integral. Dividing by a suitable power of $p$ we assume that
$$\min(v_p(u_1),v_p(u_2),v_p(u_3),v_p(c))=0.$$ We need further reductions,
as follows.

\begin{lemma}\label{lemred} Assume as above that
$\min(v_p(u_1),v_p(u_2),v_p(u_3),v_p(c))=0$.
\begin{enumerate}\item If $\min(v_p(u_1),v_p(u_2),v_p(u_3))>0$ the equation is
soluble in $\Q_p$.
\item Assume that $v_p(c)>0$, so that $\min(v_p(u_1),v_p(u_2),v_p(u_3))=0$.
The equation is equivalent to one where either $v_p(c)=0$ or
$$\max(v_p(u_1),v_p(u_2),v_p(u_3))\le2\;.$$
\item Assume that $v_p(c)>0$ and that $\max(v_p(u_1),v_p(u_2),v_p(u_3))\le2$
(which can be achieved by (2)), and without loss of generality order the
variables so that $0=v_p(u_1)\le v_p(u_2)\le v_p(u_3)\le2$, and let
$\bv=(v_p(u_2),v_p(u_3))$. Then
\begin{enumerate}\item If $\bv=(1,2)$ the equation is not soluble in $\Q_p$.
\item Or the equation is equivalent to one such that either $v_p(c)=0$
or $v_p(u_1u_2u_3)\le2$ (and $\min(v_p(u_1),v_p(u_2),v_p(u_3))=0$ if
$v_p(c)>0$).\end{enumerate}\end{enumerate}\end{lemma}

Thus, given a completely general cubic equation of the form 
$u_1X^3+u_2Y^3+u_3Z^3-cXYZ=0$, we use the following procedure, where we
distinguish the cases $c\ne0$ and $c=0$.

\smallskip

Assume first that $c\ne0$.

\begin{enumerate}\item Let $g=\gcd(u_1,u_2,u_3,c)$, and replace 
$(u_1,u_2,u_3,c)$ by $(u_1/g,u_2/g,u_3/g,c/g)$, so that we may now assume that
$\gcd(u_1,u_2,u_3,c)=1$.
\item For each prime $p\mid c$, do the following.
\begin{enumerate}
\item By dividing $u_1$, $u_2$, $u_3$, and $c$ by suitable powers of $p$
as explained in the lemma, we reduce to an equation with either
$v_p(c)=0$ or $\max(v_p(u_1),v_p(u_2),v_p(u_3))\le2$. 
\item If now $v_p(c)=0$ or if $\min(v_p(u_1),v_p(u_2),v_p(u_3))>0$ we do 
nothing more for the prime $p$. Otherwise, reorder the variables $u_i$ so that
$0=v_p(u_1)\le v_p(u_2)\le v_p(u_3)\le2$.
\item If $v_p(u_2)=1$ and $v_p(u_3)=2$, the equation has no solution.
\item Otherwise, if necessary by changing $(u_1,u_2,u_3,c)$ into
$$(u_2/p^2,u_3/p^2,pu_1,c/p)\;,$$ we may assume that we also have
$v_p(u_1u_2u_3)\le 2$.
\end{enumerate}\end{enumerate}

\smallskip

Assume now that $c=0$.

\begin{enumerate}\item Let $g=\gcd(u_1,u_2,u_3)$, and replace 
$(u_1,u_2,u_3)$ by $(u_1/g,u_2/g,u_3/g)$, so that we may now assume that
$\gcd(u_1,u_2,u_3)=1$.
\item Replace $u_1$, $u_2$, and $u_3$ by their cubefree part, so that
$\max(v_p(u_1),v_p(u_2),v_p(u_3))\le2$.
\item For each prime $p\mid u_1u_2u_3$, do the following.
\begin{enumerate}\item Reorder the variables $u_i$ so that
$0=v_p(u_1)\le v_p(u_2)\le v_p(u_3)\le2$.
\item If $v_p(u_2)=1$ and $v_p(u_3)=2$, the equation has no solution.
\item Otherwise, if necessary by changing $(u_1,u_2,u_3)$ into
$$(u_2/p^2,u_3/p^2,pu_1)\;,$$ we may assume that we also have
$v_p(u_1u_2u_3)\le 2$.
\end{enumerate}\end{enumerate}

\smallskip

This leads to the following definition:

\begin{definition} We will say that a cubic equation is \emph{$p$-reduced}
if $$\min(v_p(u_1),v_p(u_2),v_p(u_3))=0$$ and if for all primes $p$ dividing
$c$ (all primes if $c=0$) we have $v_p(u_1u_2u_3)\le2$.\end{definition}

Thanks to the above lemma, we can therefore always assume that our cubic is
$p$-reduced, since if $\min(v_p(u_1),v_p(u_2),v_p(u_3))>0$ the cubic has
a $p$-adic solution.

\begin{lemma}\label{lemsing} Let $p$ be a prime and let 
$u_1X^3+u_2Y^3+u_3Z^3-cXYZ=0$ be a cubic with $p$-integral coefficients,
not necessarily $p$-reduced. If $p\ne3$, $v_p(u_1)=v_p(u_2)=v_p(u_3)=0$,
and $v_p(27u_1u_2u_3-c^3)=0$, the cubic is soluble in $\Q_p$.
\end{lemma}

\Proof Let us look at the singular points of the cubic over $\F_p$. First, a 
point with $Z=0$ is singular if $u_1X^3+u_2Y^3=0$, $3u_1X^2=0$, and 
$3u_2Y^2=0$, and since we assume $p\ne3$ and $v_p(u_i)=0$,
this implies $X=Y=0$, which is not possible. Thus, any singular point has 
$Z\ne0$, so we may assume that $Z=1$. Since $p\ne3$, the equation has a 
singular point in $\F_p$ for $Z=1$ if and only if $3u_1X^2-cY=0$, 
$3u_2Y^2-cX=0$, and $3u_3-cXY=0$. If there is such a singular point we 
cannot have $c=0$, otherwise $u_3=0$, in other words $v_p(u_3)\ge1$, a
contradiction. Thus $Y=3u_1X^2/c$, $X=3u_2Y^2/c=27u_1^2u_2X^4/c^3$,
hence either $X=0$, which is not possible since otherwise $X=Y=0$ hence 
$u_3=0$, or $X^3=c^3/(27u_1^2u_2)$, so that
$3u_3=cXY=3u_1X^3=c^3/(9u_1u_2)$, in other words $27u_1u_2u_3-c^3=0$, which
is also excluded. Thus the cubic is nonsingular over $\F_p$. Since it is a 
curve of genus $1$ and $3-2\sqrt{2}>0$, it follows from the Weil bounds
that for every prime $p$ it has a nontrivial point in $\F_p$.
If $p$ is not in the excluded list this point is necessary nonsingular,
and since we assume $p\ne3$ we can perform a Hensel lift to $\Z_p$ as soon
as we know that there is a solution modulo $p$, proving the lemma.\fp

\subsection{Local Solubility for $p\mid u_1u_2u_3$, $p\ne3$}

Thanks to Lemma \ref{lemred}, we may assume that our cubic is $p$-reduced,
and thanks to Lemma \ref{lemsing}, it is enough to consider the primes
$p$ such that $v_p(u_1u_2u_3)>0$, $v_p(27u_1u_2u_3-c^3)>0$, or $p=3$.
We begin by primes $p\ne3$ such that $v_p(u_1u_2u_3)>0$.
For such primes, by symmetry we may assume that $v_p(u_1)>0$, and since
the cubic is $p$-reduced we have $\min(v_p(u_2),v_p(u_3))=0$, so again by
symmetry we may assume that $0\le v_p(u_1)\le v_p(u_2)\le v_p(u_3)$.

\begin{lemma}\label{lemu1u2u3} Let $p$ be a prime, assume that our cubic
is $p$-reduced, and assume that $p\ne3$ and $v_p(u_1u_2u_3)>0$, with
$0\le v_p(u_1)\le v_p(u_2)\le v_p(u_3)$. The cubic is soluble in 
$\Q_p$ if and only if one of the following conditions is satisfied.
\begin{enumerate}\item $v_p(c)=0$.
\item $v_p(c)>0$, $v_p(u_1)=v_p(u_2)=0$, and the class of $u_1/u_2$ modulo $p$ 
is a cube in $\F_p^*$.
\item $v_p(c)>0$, $v_p(u_1)=0$, $v_p(u_2)=v_p(u_3)=1$, and the class of
$u_2/u_3$ modulo $p$ is a cube in $\F_p^*$.
\end{enumerate}\end{lemma}

\smallskip

\noindent
{\bf Remarks.}\begin{enumerate}\item Note that since the cubic is $p$-reduced,
the above lemma covers all possible cases for which $p\ne3$ and
$v_p(u_1u_2u_3)>0$: indeed, if $v_p(c)>0$ we have necessarily 
$v_p(u_1u_2u_3)\le2$, so up to ordering we have either $v_p(u_1)=v_p(u_2)=0$
(and $v_p(u_3)\le2$), or $v_p(u_1)=0$ and $v_p(u_2)=v_p(u_3)=1$.
\item It follows from the proof that the cubic is also soluble in case (1)
when $p=3$, in other words if $v_3(u_1u_2u_3)>0$ and $v_3(c)=0$, but the 
assumption $p\ne3$ is necessary in cases (2) and (3).
\end{enumerate}

\smallskip

\subsection{Local Solubility for $p\mid(27u_1u_2u_3-c^3)$, $p\ne3$}

In this section, we assume that $p$ is a prime different from $3$ such that
$p\mid(27u_1u_2u_3-c^3)$. We may also assume that $p\nmid u_1u_2u_3$ since 
these primes have already been taken care of in the preceding subsection.

\begin{lemma} Let $p$ be a prime, assume that our cubic is $p$-reduced, and 
assume that $p\ne3$, $v_p(u_1u_2u_3)=0$, and $v_p(27u_1u_2u_3-c^3)>0$.
The cubic is soluble in $\Q_p$ if and only if $u_2/u_1$ is
a cube in $\F_p^*$.\end{lemma}

\subsection{Local Solubility for $p=3$}

Finally we consider local solubility at the prime $p=3$. By the remarks made
above, when $v_3(c)=0$ we have seen that the cubic is locally soluble at $3$
if $v_3(u_1u_2u_3)>0$. We may therefore assume either that $v_3(c)>0$, or that
$v_3(c)=v_3(u_1u_2u_3)=0$. In the latter case the result is immediate:

\begin{lemma} If $v_3(c)=0$ the cubic has a solution in $\Q_3$.\end{lemma} 

% \Proof If $v_3(u_1u_2u_3)>0$ we have seen that the result is true, so assume
% that $v_3(u_1u_2u_3)=0$. It is clear that the congruence 
% $u_1X^3+u_2Y^3+u_3Z^3-cXYZ\equiv0\pmod{3}$ has
% the solution $(X_0,Y_0,Z_0)=(-u_2,u_1,0)$, and
% the partial derivative with respect to $Z$ is $3u_2Z^2-cXY$, so at
% the point $(X_0,Y_0,Z_0)$ it is equal to $cu_1u_2$, hence it is not congruent
% to $0$ modulo $3$ since we have assumed $v_3(c)=v_3(u_1u_2)=0$, so we conclude
% by Hensel's lemma.\fp

The final case to be treated is the case $v_3(c)>0$. In this case, we need
a small strengthening of Hensel's lemma, which we give in a slightly more
general form that we will need below.

\begin{lemma}\label{hen3} Set $P_0=(X_0,Y_0,Z_0)$, and let $k\ge1$. Assume
that $v_3(F(P_0))\ge 2k$ and that
$\min(v_3(F'_X(P_0)),v_3(F'_Y(P_0)),v_3(F'_Z(P_0)))=k$. Assume that all second and third partial derivatives of $F$ are divisible by $3$ at 
the point $P_0$, the condition on the third derivatives being required only
if $k=1$. There exists a $3$-adic point $P=(X,Y,Z)$ such
that $F(P)=0$ with $P\equiv P_0\pmod{3^k}$.\end{lemma}

Now for $F(P)=u_1X^3+u_2Y^3+u_3Z^3-cXYZ$ we have for instance
$F'_X(P)=3u_1X^2-cYZ$, and since $v_3(c)>0$ all the partial derivatives
are divisible by $3$ at any point, so to apply the lemma it is
enough to find a point such that $F(P_0)\equiv0\pmod{3^{2k}}$ and :
$$\min(v_3(F'_X(P_0)),v_3(F'_Y(P_0)),v_3(F'_Z(P_0)))=k.$$ We will mainly use
this lemma with $k=1$, but we will need it also with $k=2$.

In fact we need a variation of the above lemma for $k=2$.

\begin{lemma}\label{hen3str} Let $P_0=(X_0,Y_0,Z_0)$ be such that 
$$v_3(F(P_0))\ge 3\text{\quad and\quad} \min(v_3(F'_X(P_0)),v_3(F'_Y(P_0)),v_3(F'_Z(P_0)))=2\;,$$
and assume that all second, third, and fourth partial derivatives of $F$ are 
divisible by $3$ at the point $P_0$. Assume in addition that for all
$P_1\equiv P_0\pmod{3}$ such that $v_3(F(P_1))\ge3$ we also have
$\min(v_3(F'_X(P_1)),v_3(F'_Y(P_1)),v_3(F'_Z(P_1)))=2$.
There exists a $3$-adic point $P=(X,Y,Z)$ such that $F(P)=0$ with 
$P\equiv P_0\pmod{3}$.\end{lemma}

Of course for our cubics the fourth partial derivatives vanish.

\smallskip

% \Proof Set $P_1=P_0+3(u,v,w)$ for some $3$-adic integers $u$, $v$, $w$. Once
% again, by Taylor's formula, we have
% \begin{align*}F(P_1)&=F(P_0)+3(u,v,w)\cdot(F'_X(P_0),F'_Y(P_0),F'_Z(P_0))\\
% &\phantom{=}+(9/2)S_2+(27/6)S_3+\sum_{n\ge4}(3^n/n!)S_n\;,\end{align*}
% where $S_n$ is as above. We check that $v_3(3^{n-4}/n!)\ge0$ for $n\ge5$,
% so $\sum_{n\ge5}(3^n/n!)S_n\equiv0\pmod{81}$, and since $3\mid S_4$ we have
% $(3^4/4!)S_4\equiv0\pmod{81}$. Thus the congruence $F(P_1)\equiv0\pmod{81}$
% is equivalent to 
% $$(u,v,w)\cdot(F'_X(P_0)/9,F'_Y(P_0)/9,F'_Z(P_0)/9)\equiv 
% -(F(P_0)/27+S_2/6+S_3/6)\pmod{3}\;,$$
% which is possible since by assumption one of the first partial derivatives
% has valuation equal to $2$. Thus we can solve the congruence 
% $F(P_1)\equiv0\pmod{81}$. By the additional assumption that we have made,
% we know that $$\min(v_3(F'_X(P_1)),v_3(F'_Y(P_1)),v_3(F'_Z(P_1)))=2\;,$$ so
% we may apply Lemma \ref{hen3} with $k=2$ to the point $P_1$, which tells us
% that there exists a $3$-adic solution congruent to $P_1$, hence also
% to $P_0$, modulo $3$.\fp

Recall that a solution $(X,Y,Z)$ of a congruence modulo some power of $3$ 
is always such that $\min(v_3(X),v_3(Y),v_3(Z))=0$.

\begin{lemma}\label{lemp3} Let $p=3$, assume the cubic $3$-reduced, and assume
that $v_3(c)>0$, so that $v_3(u_1u_2u_3)\le2$. Reorder the variables so that 
$0=v_3(u_1)\le v_3(u_2)\le v_3(u_3)$.
\begin{enumerate}\item If $v_3(c)\ge2$ and $v_3(u_1u_2u_3)=0$ the cubic has a 
solution in $\Q_3$ if and only if $u_i\equiv \pm u_j\pmod{9}$ for some 
$i\ne j$.
\item If $v_3(c)\ge2$ and exactly one of the $u_i$ is divisible by 
$3$ (in other words if $v_3(u_2)=0$ and $v_3(u_3)>0$), the cubic has a 
solution in $\Q_3$ if and only if either $u_1\equiv\pm u_2\pmod{9}$, or if
$v_3(u_3)=1$.
\item If $v_3(c)\ge2$, and two of the $u_i$ are divisible by $3$ (in other
words if $v_3(u_2)=v_3(u_3)=1$ since the cubic is $3$-reduced), the cubic has 
a solution in $\Q_3$ if and only if $u_2/3\equiv\pm u_3/3\pmod{9}$.
\item If $v_3(c)=1$ and exactly one of the $u_i$ is divisible by $3$
(i.e., $v_3(u_2)=0$ and $v_3(u_3)>0$), the cubic has a solution in 
$\Q_3$ if and only if either $u_1\equiv \pm u_2\pmod{9}$, or if there exist 
$s_1$ and $s_2$ in $\{-1,1\}$ such that 
$c\equiv s_1u_1+s_2u_2+s_1s_2u_3\pmod{9}$.
\item If $v_3(c)=1$ and two of the $u_i$ are divisible by $3$
(i.e., $v_3(u_2)=v_3(u_3)=1$), the cubic has a solution in $\Q_3$.
\end{enumerate}\end{lemma}

Note that this lemma does not cover the case where $v_3(c)=1$ and 
none of the $u_i$ is divisible by $3$, or equivalently $v_3(u_1u_2u_3)=0$,
which will be covered by Lemma \ref{a3bnot} below.

\smallskip

\begin{lemma}\label{a3bnot} Let $p=3$ and assume that $v_3(c)>0$ and
$v_3(u_1u_2u_3)=0$.
\begin{enumerate}\item If $u_i\equiv \pm u_j\pmod{9}$ for some $i\ne j$, the 
cubic has a solution in $\Q_3$.
\item If $u_i\not\equiv \pm u_j\pmod{9}$ for $i\ne j$ (which implies that
$u_1u_2u_3\equiv\pm1\pmod{9}$), the cubic has a solution in 
$\Q_3$ if and only if there exist suitable signs $s_1=\pm1$ and $s_2=\pm1$
such that $c\equiv s_1u_1+s_2u_2+s_1s_2u_3\pmod{27}$.
\end{enumerate}\end{lemma}

\smallskip

\noindent
{\bf Remarks.}\begin{enumerate}\item We only assume that $v_3(c)>0$ and not
$v_3(c)=1$, although the case $v_3(c)\ge2$ is covered by Lemma \ref{lemp3}:
indeed, it is easy to see case (2) of the above lemma cannot occur when
$v_3(c)\ge2$.
\item It follows from the proof that case (1) occurs if and only if there
exists a solution $(X,Y,Z)$ with $\min(v_3(X),v_3(Y),v_3(Y))=0$ but one
of the variables divisible by $3$.
\item In case (2), there is no need to search among the four possibilities
for the signs $s_i$: since $3\mid c$ it is easy to see that we must take
$s_1\equiv u_2u_3\pmod{3}$ and $s_2\equiv u_1u_3\pmod{3}$.
\end{enumerate}

\smallskip

\smallskip

We have thus finished to give the local solubility of the general
cubic equation $$F(X,Y,Z)=u_1X^3+u_2Y^3+u_3Z^3-cXYZ=0,$$ hence in particular
of the equation $$u_1X^3+u_2Y^3+(2b/(u_1u_2))Z^3-2aXYZ=0$$ of Theorem
\ref{thd1}.

%%%%%%%%%%%%%%%%%%%%%%%%%%%%%%%%%%%%%%%%%%%%%%%%%%%%%%%%%%%%%%%%%%%%%%%

\section{Local Solubility : The Case $D\ne1$}

\subsection{Reduction of the Cubic, Bad Primes, and Split Primes}

We now consider local solubility when $D\ne1$. Although we will do the same
type of computations as in the case $D=1$, there are evidently some added 
complications.

It is essential to begin by a reduction of the cubic equation of Theorem
\ref{thdnot1}. Recall that in the case $D=1$ we could reduce to an equation
where $u=u_1^2u_2$ with the $u_i$ squarefree and coprime. We have seen
in Section \ref{secg3} that the analogous statement for $D\ne1$ involves
ideals, so we cannot immediately reproduce what we have done.

Recall that the cubic equation of Theorem \ref{thdnot1} can be written as
$F(X,Y,Z)=0$, where if $v=v_1+v_2\sqrt{D}$ we have
\begin{align*}F(X,Y,Z)&=2v_2X^3+2Dv_1Y^3+(2b/(v_1^2-Dv_2^2))Z^3\\
&\phantom{=}+6v_1X^2Y+6v_2DXY^2+2a(X^2Z-DY^2Z)\\
&=(v(X+Y\sqrt{D})^3-\tau(v)(X-Y\sqrt{D})^3)/\sqrt{D}\\
&\phantom{=}+2aZ(X+Y\sqrt{D})(X-Y\sqrt{D})+(2b/(v\tau(v)))Z^3\;,\end{align*}
and we will use indifferently both forms.

\begin{lemma}\label{lemsingnot1} If $p\ne3$, $v_p(v\tau(v))=0$, $v_p(2b)=0$,
and $v_p(27b-4a^3D)=0$ the above cubic is soluble in $\Q_p$.
\end{lemma}

% \Proof Since the above cubic comes from that of Theorem \ref{thd1} after some
% simple transformations which do not add any bad primes except possibly $p=2$,
% this is a direct consequence of Lemma \ref{lemsing}. More precisely, the 
% simplest way to go about the proof is to set $\al=Y+X\sqrt{D}$ and 
% $\be=Y-X\sqrt{D}$, considered as independent variables. The (affine) equation
% is then $(v\al^3-\tau(v)\be^3)/\sqrt{D}+2a\al\be+2b/(v\tau(v))=0$, and
% computing derivatives with respect to $\al$ and $\be$ as we have done for Lemma
% \ref{lemsing} gives the desired result. Note that that lemma is a special 
% case of the present one if we understand $v\tau(v)$ as $u_1u_2$, which is 
% reasonable since here $u=v^2\tau(v)$, while $u=u_1^2u_2$ in Lemma
% \ref{lemsing}.\fp

\smallskip

Recall that, analogously to ideals, an algebraic integer $v$ is said to be
primitive if $v/n\in\Z_K$ with $n\in\Z$ if and only if $n=\pm1$.

\begin{lemma}\label{lemsp} Let $[u]\in\Ima(\al)$. In the above cubic, we may 
assume that $[u]=[v^2\tau(v)]$ where $v$ is a primitive algebraic integer such
that $v\tau(v)$ is divisible only by split primes. In particular, $v$ and 
$\tau(v)$ generate coprime ideals. Furthermore, if $D\equiv0\pmod{4}$ we
may also assume that $v=v_1+v_2\sqrt{D}$ with $v_2\in\Z$.
\end{lemma}

\smallskip

Note that, contrary to the case $D=1$, we cannot deduce from this lemma
that $v\tau(v)\mid 2b$. It is easy to show using Corollary \ref{corq3}
that if $3\nmid h(K)$, we may assume that $v\tau(v)\mid(2b)^{h(K)}$,
and in particular all prime numbers dividing $v\tau(v)$ (which are necessarily
split) divide $2b$, but this may not be true if $3\mid h(K)$.

Since $v$ is now an algebraic integer, we have $2v_1\in\Z$ and $2v_2\in\Z$,
so all the coefficients of the equation are integers, except perhaps for
$2b/(v\tau(v))$. 

\begin{corollary}\label{corsplit} Assume as above that $v=v_1+v_2\sqrt{D}$ is
a primitive algebraic integer such that $v\tau(v)$ is divisible only by split
primes, and let $p$ be any split prime. There exists $d_p\in\Q_p$ such that
$d_p^2=D$. The cubic of Theorem \ref{thdnot1} has a solution in $\Q_p$
if and only if the cubic $u_1X^3+u_2Y^3+u_3Z^3-cXYZ=0$ does, where
$u_1=v_1+v_2d_p$, $u_2=v_1-v_2d_p$, $u_3=(2b/v\tau(v))d_p$, and $c=2ad_p$.
\end{corollary}

% \Proof Since $p$ is split we have $\lgs{D}{p}=1$, so there exists $d$ such
% that $d^2\equiv D\pmod{p}$, and $p\nmid d$, and since $D$ is assumed to be a
% discriminant, we may assume if necessary that $d^2\equiv D\pmod{4p}$, which
% is essential for $p=2$. Using Hensel lifting, we can now lift $d$
% to a true square root of $D$ in $\Z_p$, including for $p=2$, so there exists
% $d_p\in\Z_p$ such that $d_p^2=D$. Since we are looking for local solutions, we
% may replace $\sqrt{D}$ by $d_p$ without changing the problem, and of course if
% $v=v_1+v_2\sqrt{D}$ the corresponding $p$-adic numbers are $v=v_1+v_2d_p$ and 
% $\tau(v)=v_1-v_2d_p$. Thus, if $p$ is any split prime, we will need to solve 
% in $\Q_p$ the equation
% $$(v(X+Yd_p)^3-\tau(v)(X-Yd_p)^3)/d_p+2aZ(X+Yd_p)(X-Yd_p)+(2b/(v\tau(v)))Z^3=0\;.$$
% We can thus set $X_1=X+Yd_p$, $Y_1=-(X-Yd_p)$, and use these variables instead
% (note that since we are in characteristic $0$ we can do this also for $p=2$),
% so the equation reads, after setting $u_1=v$ and $u_2=\tau(v)$ (which are 
% independent elements of $\Z_p$):
% $$u_1X_1^3+u_2Y_1^3+(2b/(v\tau(v))d_p)Z^3-2ad_pX_1Y_1Z=0\;,$$
% proving the corollary.\fp

Since we have given a complete algorithm to determine the
solubility in $\Q_p$ of an equation of the type $u_1X^3+u_2Y^3+u_3Z^3-cXYZ=0$,
this solves the problem for $D\ne1$ in the case where $p$ is a split prime. 
Thus we are left with the study of ramified and inert primes, so thanks to 
Lemma \ref{lemsp}, we may assume that $v_p(v\tau(v))=0$, so that in
particular $v_p(2b/(v\tau(v)))\ge0$.

\smallskip

Thus in the sequel we assume that $p$ is a ramified or inert prime,
that $v_p(v\tau(v))=0$, and recall that our equation is $F(X,Y,Z)=0$ with 
\begin{align*}F(X,Y,Z)&=(v(X+Y\sqrt{D})^3-\tau(v)(X-Y\sqrt{D})^3)/\sqrt{D}+u_3Z^3\\
&\phantom{=}+2aZ(X+Y\sqrt{D})(X-Y\sqrt{D})\;,\end{align*}
with $u_3=2b/(v\tau(v))$, hence such that $v_p(u_3)\ge0$.
We begin by inert primes.

\subsection{The Case of Inert Primes}

If $p$ is an inert prime, consider the field $K_p=\Q_p(\sqrt{D})$, which up
to isomorphism is the unique unramified extension of degree $2$ of $\Q_p$,
and whose residue field is $\F_{p^2}$, so that we can also consider
the class of $\sqrt{D}$ in $\F^*_{p^2}$, and $\tau$ is defined on $\F_{p^2}$
as in characteristic $0$. By abuse of notation, if $\al$ and $\be$ are
elements of $K=\Q(\sqrt{D})$ or of $K_p$, we will write $\al\equiv\be\pmod{p}$
to mean that the class of $\al$ and $\be$ in the residue field is the same.
Note that we work in $K_p$ or $K$ for practicality, but that the cubic
equation has coefficients in $\Q$, and that we also look for solutions in $\Q$.

Before studying the bad primes, we need an auxiliary lemma.

\begin{lemma}\label{lemcubfp2} Let $p\ne3$ be an inert prime. The following
conditions are equivalent:
\begin{enumerate}\item There exist $X$ and $Y$ such that 
$\tau(v)/v\equiv ((X+Y\sqrt{D})/(X-Y\sqrt{D}))^3\pmod{p}$ in the above sense.
\item The class of $\tau(v)/v$ is a cube in $\F_{p^2}^*$.
\item We have either $p\equiv1\pmod{3}$, or $p\equiv2\pmod{3}$ and
$v^{(p^2-1)/3}\equiv1\pmod{p}$.\end{enumerate}\end{lemma}

% \Proof (1) implies (2) is trivial. Let us show conversely that (2) implies (1),
% so assume that $\tau(v)/v\equiv\ga^3\pmod{p}$ for some $\ga\in K$.
% Applying $\tau$ and multiplying we obtain that 
% $(\ga\tau(\ga))^3\equiv1\pmod{p}$, so that $\ga\tau(\ga)\equiv\rho\pmod{p}$,
% where $\rho$ is some cube root of $1$ modulo $p$. This implies that in
% $\F_{p^2}$ we have $\tau(\rho)\equiv\rho\pmod{p}$. Thus 
% $(\rho\ga)\tau(\rho\ga)\equiv\rho^3\equiv1\pmod{p}$, so if we set
% $\al=1+\tau(\rho\ga)$, we have 
% $$\tau(\al)=1+\rho\ga\equiv\rho\ga\tau(\rho\ga)+\rho\ga\equiv\rho\ga(1+\tau(\rho\ga))\equiv\rho\ga\al\pmod{p}\;,$$
% so that
% $\tau(v)/v\equiv\ga^3\equiv(\rho\ga)^3\equiv(\tau(\al)/\al)^3\pmod{p}$.
% Thus, if we set $\al=X-Y\sqrt{D}$, we obtain that
% $\tau(v)/v\equiv((X+Y\sqrt{D})/(X-Y\sqrt{D}))^3\pmod{p}$, proving (1).
% 
% \smallskip
% 
% Let us show now that (2) and (3) are equivalent. Recall that for any
% $\al\in\F_{p^2}$ we have $\tau(\al)=\al^p$. Let $g$ be a generator of the
% cyclic group $\F_{p^2}^*$, and let $k$, defined uniquely modulo $p^2-1$, be
% such that $v\equiv g^k\pmod{p}$. We thus have $\tau(v)/v\equiv g^{k(p-1)}$,
% so the class of $\tau(v)/v$ is a cube in $\F_{p^2}^*$ if and only if
% $3\mid k(p-1)$, in other words if either $p\equiv1\pmod{3}$ or if $3\mid k$,
% and the condition $3\mid k$ is evidently equivalent to 
% $v^{(p^2-1)/3}\equiv1\pmod{p}$, proving the lemma.\fp

\begin{lemma}\label{lemin2b} Let $D\ne1$, and assume that the elliptic curve 
is given by an equation satisfying the conditions of Lemmas \ref{lemeq} and
\ref{lemsp}. Let $p$ be an inert prime number such that $v_p(2b)>0$, 
$v_p(v\tau(v))=0$ and $p\ne3$, and if $p=2$, assume that $v_p(2b)\le2$. The 
cubic of Theorem \ref{thdnot1} is locally soluble at $p$ if and only if one of
the following conditions is satisfied.
\begin{enumerate}\item $v_p(2a)=0$.
\item $v_p(2a)>0$ and the class of $\tau(v)/v$ modulo $p$ 
is a cube in $\F_{p^2}^*$.
\end{enumerate}
\end{lemma}

We now consider the prime $p=2$, assumed to be inert, when $v_2(2b)\ge3$.
Since the equation is $2$-reduced, note that either $v_2(a)>0$, in which case
$v_2(b)\le2$ hence $v_2(2b)=3$, or $v_2(a)=0$. Furthermore, we can write 
$v=v_1+v_2\sqrt{D}=(w_1+w_2\sqrt{D})/2$ with $w_1$ and $w_2$ in $\Z$ such that
$w_1\equiv w_2\pmod{2}$, and since $v$ is primitive, either we have 
$w_1\equiv w_2\equiv1\pmod{2}$, or $w_1$ and $w_2$ are even with 
$w_1\not\equiv w_2\pmod{4}$. 

\begin{lemma}\label{lemin2b2} Let $D\ne1$, and assume that the elliptic curve 
is given by an equation satisfying the conditions of Lemmas \ref{lemeq} and
\ref{lemsp}, assume that $p=2$ is an inert prime, in other words that
$D\equiv5\pmod{8}$, assume that $v_2(2b)\ge3$, and write $w_1=2v_1$ and
$w_2=2v_2$.
\begin{enumerate}\item If $w_1\equiv2\pmod{4}$ and $w_2\equiv0\pmod{4}$
or $w_1\equiv0\pmod{4}$ and $w_2\equiv2\pmod{4}$ the cubic has a solution in
$\Q_2$.
\item If $w_1\equiv w_2\equiv1\pmod{2}$, the cubic has a solution in $\Q_2$ if
and only if either $v_2(2b)\ge4$ or $v_2(a)>0$.
\end{enumerate}\end{lemma}

\begin{lemma} Let $D\ne1$, and assume that the elliptic curve is given by
an equation satisfying the conditions of Lemmas \ref{lemeq} and \ref{lemsp}. 
Let $p$ be an inert prime number such that $v_p(2b)=v_p(v\tau(v))=0$, 
$v_p(27b-4Da^3)>0$, and $p\ne3$. The cubic of Theorem \ref{thdnot1} is locally 
soluble at $p$ if and only if $\tau(v)/v$ is a cube in $\F_{p^2}^*$.\end{lemma}

Recall that since we assume that $p$ is not split we have $v_p(v\tau(v))=0$.

\begin{lemma} Let $D\ne1$, assume that the elliptic curve is given by
an equation satisfying the conditions of Lemmas \ref{lemeq} and \ref{lemsp},
and assume that $p=3$ is an inert prime, i.e., that $D\equiv2\pmod{3}$.
Set $u_1=2v_2$, $u_2=2v_1D$, and $u_3=2b/(v\tau(v))$.
\begin{enumerate}\item If $v_3(2a)=0$ the cubic has a solution in
$\Q_3$.
\item If $v_3(2a)\ge2$ the cubic has a solution in $\Q_3$ if and only if 
either $v_3(u_1)\ge2$, $v_3(u_2)\ge2$, $u_i\equiv\pm u_j\pmod{9}$ for some 
$i\ne j$ and a suitable sign $\pm$, or if $u_3\equiv2(\pm u_1\pm u_2)\pmod{9}$
for suitable signs $\pm$.
\item If $v_3(2a)=1$, $v_3(2b)>0$, and $v_3(u_1u_2)\ge1$, the cubic has a 
solution in $\Q_3$ if and only if either $v_3(u_1)\ge2$, $v_3(u_2)\ge2$, or
$v_3(2a+2b)=1$.
\item If $v_3(2a)=1$, $v_3(2b)>0$, and $v_3(u_1u_2)=0$, the cubic has a 
solution in $\Q_3$ if and only if either $u_1\equiv\pm u_2\pmod{9}$,
or if $2a+b\equiv \pm u_1\pm u_2\pmod{9}$ for suitable signs $\pm$.
\end{enumerate}\end{lemma}

\begin{lemma} Let $D\ne1$, assume that the elliptic curve is given by
an equation satisfying the conditions of Lemmas \ref{lemeq} and \ref{lemsp},
and assume that $p=3$ is an inert prime, i.e., that $D\equiv2\pmod{3}$.
Set $u_1=2v_2$, $u_2=2v_1D$, and $u_3=2b/(v\tau(v))$. Assume that $v_3(2a)=1$
and $v_3(u_1)=v_3(u_2)=v_3(2b)=0$, and set $u_4=u_3+2a$, so that we also
have $v_3(u_4)=0$.
\begin{enumerate}\item If $u_i\equiv\pm u_j\pmod{9}$ for some $i\ne j$
with $i,j=1$, $2$, or $4$, the cubic has a solution in $\Q_3$.
\item Otherwise, the cubic has a solution in $\Q_3$ if and only if
$2a+u_4=4a+u_3\equiv \pm 3D\pmod{27}$ for a suitable sign $\pm$.
\end{enumerate}\end{lemma}

Similar remarks as those given after Lemma \ref{a3bnot} apply here.

\smallskip

This concludes the study of local solubility in the case of inert primes.

\subsection{The Case of Ramified Primes}

\begin{lemma} Let $D\ne1$, and assume that the elliptic curve is given by
an equation satisfying the conditions of Lemmas \ref{lemeq} and \ref{lemsp}. 
If $p$ is a ramified prime such that $p\ne3$, the cubic is soluble in $\Q_p$.
\end{lemma}

% \Proof Recall that 
% \begin{align*}F(X,Y,Z)&=2v_2X^3+2Dv_1Y^3+(2b/(v_1^2-Dv_2^2))Z^3\\
% &\phantom{=}+6v_1X^2Y+6v_2DXY^2+2a(X^2Z-DY^2Z)\;.\end{align*}
% Assume first that $p\ne2$. If $p$ is ramified we have $p\mid D$, and since 
% $v\tau(v)$ is only divisible by split primes we have $v_p(v_1^2-Dv_2^2)=0$, 
% and in particular $v_p(2v_1)=0$ since $p\ne2$. Thus
% $F(X,Y,Z)\equiv 2v_2X^3+(2b/(v\tau(v)))Z^3+6v_1X^2Y+2aX^2Z\pmod{p}$.
% Since $v_p(2v_1)=0$ and $p\ne3$, we choose $Z_0=0$, $X_0=1$, and
% $Y_0\equiv-2v_2/(6v_1)\pmod{p}$, so that $F(X_0,Y_0,Z_0)\equiv0\pmod{p}$. 
% Furthermore, $F'_Y(X_0,Y_0,Z_0)\equiv 6v_1\not\equiv0\pmod{p}$, so
% we conclude by Hensel's lemma that the cubic is soluble in $\Q_p$.
% 
% Assume now that $p=2$. Since $p$ is ramified we have $D\equiv0\pmod{4}$,
% hence by Lemma \ref{lemsp} we may assume that $v_1$ and $v_2$ are in $\Z$.
% Since $v\tau(v)=v_1^2-Dv_2^2$ is odd and $4\mid D$, it follows that $v_1$
% is odd. Thus we must solve $G(X,Y,Z)=0$, where
% \begin{align*}G(X,Y,Z)&=F(X,Y,Z)/2=v_2X^3+Dv_1Y^3+(b/(v_1^2-Dv_2^2))Z^3\\
% &\phantom{=}+3v_1X^2Y+3v_2DXY^2+a(X^2Z-DY^2Z)\end{align*}
% has $2$-integral coefficients. We reason as for $p\ne2$: choosing $Z_0=0$,
% $X_0=1$, and $Y_0=-v_2/(3v_1)\pmod{2}$, we have $G(X_0,Y_0,Z_0)\equiv0\pmod{2}$
% and $G'_Y(X_0,Y_0,Z_0)\equiv 3v_1\not\equiv0\pmod{2}$, so once again we 
% conclude by Hensel's lemma.\fp

\begin{lemma}\label{lem3ram1} Let $D\ne1$, and assume that the elliptic curve
is given by an equation satisfying the conditions of Lemmas \ref{lemeq} and 
\ref{lemsp}. Assume that $p=3$ is ramified, in other words that $3\mid D$, and
to simplify notation, set $u_3=2b/(v\tau(v))$.
\begin{enumerate}\item If $v_3(2a)=0$, the cubic has a solution
if and only if one of the following conditions is satisfied:
\begin{enumerate}\item $v_3(2v_2)>0$.
\item $v_3(2v_2)=v_3(2a+u_3)=0$.
\end{enumerate}
\item If $v_3(2a)\ge2$, the cubic has a solution if and only one of
the following conditions is satisfied:
\begin{enumerate}
\item $D\equiv 3\pmod{9}$ and $v_3(u_3)=0$.
\item $D\equiv 3\pmod{9}$, $v_3(u_3)>0$, and $v_3(2v_2)>0$.
\item $D\equiv 6\pmod{9}$ and $v_3(2v_2)\ge2$.
\item $D\equiv 6\pmod{9}$ and $v_3(2v_2)=v_3(u_3)=1$.
\item $D\equiv 6\pmod{9}$, $v_3(2v_2)=0$, and $u_3\equiv\pm 2v_2\pmod{9}$.
\item $u_3\equiv\pm 2v_1D\pmod{27}$.
\end{enumerate}\end{enumerate}\end{lemma}

We will treat the case where $v_3(2a)=1$ below.

\smallskip

\begin{lemma} Keep the notation and assumptions of the preceding lemma,
assume now that $v_3(2a)=1$, and set $u_4=u_3+2a$. The cubic has a solution if
and only one of the following conditions is satisfied:
\begin{enumerate}
\item[(a)] $D\equiv 3\pmod{9}$ and $v_3(u_4)=0$.
\item[(b)] $D\equiv 3\pmod{9}$, $v_3(u_4)>0$, and $v_3(2v_2)>0$.
\item[(c)] $D\equiv 6\pmod{9}$ and $v_3(2v_2)\ge2$.
\item[(d)] $D\equiv 6\pmod{9}$ and $v_3(2v_2)=v_3(u_4)=1$.
\item[(e)] $D\equiv 6\pmod{9}$, $v_3(2v_2)=0$, and $u_4\equiv\pm 2v_2\pmod{9}$.
\item[(f)] $v_3(u_3)=1$ and there exists $s=\pm1$ such that
$2v_1(D/3)\equiv s(u_3/3-2a(D/3))\pmod{9}$ and $2v_1s\not\equiv 2a/3\pmod{3}$.
\item[(g)] $v_3(u_3)=1$ and there exists $s=\pm1$ such that
$2v_1(D/3)\equiv s(u_3/3-2a(D/3))\pmod{9}$, $2v_1s\equiv 2a/3\pmod{3}$,
$v_3(2v_2)=0$, and $D\equiv 3\pmod{9}$.
\item[(h)] $v_3(u_3)=1$ and there exists $s=\pm1$ such that
$2v_1(D/3)\equiv s(u_3/3-2a(D/3))\pmod{27}$, $2v_1s\equiv 2a/3\pmod{3}$,
$v_3(2v_2)=0$, and $D\equiv 6\pmod{9}$.
\item[(i)] $v_3(u_3)=1$ and there exists $s=\pm1$ such that
$2v_1s\equiv 2a/3\pmod{3}$, $v_3(2v_2)>0$, and there exists $t\in\{-1,0,1\}$ 
and $r\in\{-1,0,1\}$ such that
\begin{align*}2v_1(D/3)&\equiv s(u_3/3-2a(D/3))-6v_2(D+3)t-9(2v_1s-2a/3)st^2\\
&\phantom{=}-3r(D(2v_1+a/3)+6rDv_1+6at^2)\pmod{81}\;.\end{align*}
\end{enumerate}
\end{lemma}

\smallskip

This finishes the study of local solubility in the case of ramified primes,
hence the section on local solubility.

\section{Examples}

\subsection{The Curves $y^2=x^3+(kp)^2$ for $k=1$, $2$, or $4$}

In this section we consider the family of curves $E_{kp}$ with equation
$y^2=x^3+(kp)^2$, where $p$ is a prime and $k=1$, $2$, or $4$. The restriction
on $k$ is made so that no other prime apart from $2$ divides it. Note that
it is not necessary to consider higher powers of $2$ since the
curve $y^2=x^3+(8kp)^2$ is trivially isomorphic to the curve $y^2=x^3+(kp)^2$.
Furthermore the primes $p=2$ and $3$ give rise to a finite number of curves
which can be treated individually (specifically, for $p=2$ the rank is equal
to $0$, for $(k,p)=(1,3)$ and $(2,3)$ the rank is $1$, Mordell--Weil generators
being $(-2,1)$ and $(-3,3)$ respectively, and for $(k,p)=(4,3)$ the rank is
again equal to $0$, and the torsion is always of order $3$ generated by
$T=(0,kp)$, except for $(k,p)=(4,2)$, for which it has order $6$ generated by 
$(8,24)$). We therefore assume that $p\ge5$, so that in particular all of
these curves have rational $3$-torsion generated by $T=(0,kp)$ equal to
their full rational torsion subgroup.

\smallskip

We first compute the image of $\al$. For this, we consider the cubic equations
of Theorem \ref{thd1} (3), in other words $u_1X^3+u_2Y^3+u_3Z^3=0$,
where $u_1u_2\mid 2kp$ and $u_3=2kp/(u_1u_2)$, where we recall that $u_1u_2$
is squarefree. Up to exchange of $u_1$ and $u_2$, it is easy to check that
the only of possibilities are $(1,1,2kp)$, $(1,2,kp)$, $(1,p,2k)$, $(1,2p,k)$,
and $(2,p,k)$. The first one (corresponding to $u=1$) gives an evidently
soluble equation, corresponding to the unit element of the elliptic curve.

\begin{itemize}\item When $k=1$, the fourth one (corresponding to $u=2p$ and 
$u=4p^2$) is also soluble, and it corresponds to the two nontrivial rational
$3$-torsion points on the curve, and the other three (corresponding to
$u=2$, $4$, $p$, $p^2$, $4p$, and $2p^2$, are equivalent.
\item When $k=2$, the fifth one (corresponding to $u=4p$ and $u=2p^2$) is
also soluble, and it again corresponds to the two nontrivial rational
$3$-torsion points on the curve, the second and fourth (corresponding to
$u=2$, $4$, $2p$, and $4p^2$ are equivalent, and the third corresponds to
$u=p$ and $p^2$. Since the set (of classes) of $u$ for which the cubic is 
soluble forms a group and since $4p$ and $2p^2$ belong to this group, 
it follows that $[p]$ and $[p^2]$ will be in the group if and only if
$4$ and $2$ are, so in fact the three equations are equivalent, although
slightly less trivially.
\item Finally, when $k=4$ the third equation (corresponding to $p$ and $p^2$)
is clearly soluble (since $2k=8$ is a cube), and this again corresponds to
the two rational $3$-torsion points. The other equations correspond 
respectively to $u=2$ and $4$, $u=2p$ and $4p^2$, and $u=4p$ and $2p^2$, and
once again because of the group structure all the equations are in fact
equivalent. 
\end{itemize}

We see that in each case it is sufficient to consider the equation with
$u_1=1$, $u_2=2$, hence $u_3=kp$. The result is as follows:

\begin{lemma} Keep the above assumptions. The equation $X^3+2Y^3+kpZ^3=0$ is
ELS if and only if $k\ne4$, either $p\equiv2\pmod{3}$ or 
$2^{(p-1)/3}\equiv1\pmod{3}$, and $kp\not\equiv\pm4\pmod{9}$.\end{lemma}

\Proof This of course immediately follows from the study of the equation that
we have made above. More precisely, if $k=4$ the $2$-adic valuations of the
coefficients are $(0,1,2)$, so the equation has no $2$-adic solutions by
Lemma \ref{lemred}. On the other hand, if $k=1$ or $k=2$, the $2$-adic
valuations are $(0,0,1)$ and $(0,1,1)$ respectively, and since all elements
of $\F_2^*$ are cubes, we conclude by Lemma \ref{lemu1u2u3} that the equation
has a solution in $\Q_2$. For $p$-adic solubility we also use this lemma,
since the $p$-adic valuations are $(0,0,1)$, and we conclude that the equation
has a $p$-adic solution if and only if $2$ is a cube in $\F_p^*$, leading to
the given condition. Finally, since the $3$-adic valuations are $(0,0,0)$
we use Lemma \ref{lemp3} (1), which tells us that the equation has a $3$-adic
solution if and only if $u_i\equiv\pm u_j\pmod{9}$ for some $i\ne j$,
which gives $kp\equiv\pm1$ or $\pm2\pmod{9}$, in other words 
$kp\not\equiv\pm4\pmod{9}$ since $3\nmid kp$.\fp

\begin{corollary}\label{corex1} Let $p\ge5$ be prime, let $k=1$, $2$, or $4$,
and let $E$ be the elliptic curve $y^2=x^3+(kp)^2$.
\begin{enumerate}\item For $k=1$, if either $p\equiv\pm4\pmod{9}$ or if 
$p\equiv1$ or $7\pmod9$ and $2^{(p-1)/3}\not\equiv1\pmod{p}$ then
$\Ima(\al)=\{1,2p,4p^2\}$, and in particular $|\Ima(\al)|=3$.
\item For $k=1$ and $p\equiv2\pmod{9}$ we have $|\Ima(\al)|=9$.
\item For $k=2$, if either $p\equiv\pm2\pmod{9}$, or if $p\equiv 1$ or
$4\pmod{9}$ and $2^{(p-1)/3}\not\equiv1\pmod{p}$ then
$\Ima(\al)=\{1,4p,2p^2\}$, and in particular $|\Ima(\al)|=3$.
\item For $k=4$ we always have $\Ima(\al)=\{1,p,p^2\}$, and in particular 
$|\Ima(\al)|=3$.
\item In all other cases we have $|\Ima(\al)|=3$ or $9$. More precisely
the cubic equation $X^3+2Y^3+kpZ^3=0$ is ELS, and
$|\Ima(\al)|=9$ if and only if it is globally soluble.
\end{enumerate}
\end{corollary}

\Proof (1), (3), (4), and (5) are clear from the lemma by inspection. For (2),
we use the proposition 3.3 p. 438 of Satg\'e \cite{Satg2}. \fp
%Using a Heegner point construction. Here is the argument : the proposition 1.16 p. 431 from \cite{Satg2} provides the equality $f_{1}(p\rho)^3=f_{2}(p\rho)^3+2\;\;(*)$. The numbers $\alpha=f_{1}(p\rho)/\sqrt[3]{4}$ and $\beta=f_{1}(p\rho)/\sqrt[3]{2p}$ are algebraic over $\mathbb{Q}$ and Satgé's theorem 2.1 (p. 433) shows that $[\mathbb{Q}[\alpha,\beta]:\mathbb{Q}]$ is prime to $3$. If we write $(*)$ in the form $2(-f_{1}(p\rho)/\sqrt[3]{4})^3+1=p(-f_{2}(p\rho)/\sqrt[3]{2p})^3$, then we can say that the projective equation $2Y^3+X^3=pZ^3$ has a rational solution on the field $\mathbb{Q}[\alpha,\beta]$. It also has a rational solution on the field $\mathbb{Q}[\sqrt[3]{p}]$. Then it suffices to apply the lemma 3.2 p. 438 of \cite{Satg2} to conclude that it has a $\mathbb{Q}$-rational solution. It implies that the cubic $X^3+2Y^3+pZ^3=0$ is globally soluble when $p\equiv2\pmod9$, so (2) follows from (5).

\smallskip

\noindent
{\bf Remark.} In \cite{RoZa}, Rodriguez-Villegas and Zagier have characterized
the primes which are sums of two cubes. If their method could be extended to
primes which are of the form $x^3+2y^3$, and also of the form $x^3+4y^3$,
it would determine $\Ima(\al)$ in all cases.

\medskip

We now compute the image on the dual curve $\wh{E}$, whose equation is
$y^2=x^3-27(kp)^2$, so that $D=-3$ and $b=3kp$. We first determine local
solubility of the equation corresponding to $\rho=(-1+\sqrt{-3})/2$,
and for the moment we do not necessarily assume that $k\mid 4$.

\begin{lemma}\label{lemrholocold}
\begin{enumerate}\item Let $k$ be such that $8\nmid k$. The equation
corresponding to $\rho$ is locally soluble at the primes $2$, $3$, and 
$p$ if and only if $p\equiv\pm1\pmod{9}$, $k\equiv\pm4\pmod{9}$, and 
$4\mid k$.
\item In particular, if $k=1$, $2$, or $4$, the equation corresponding to
$\rho$ is ELS if and only if $k=4$ and $p\equiv\pm1\pmod{9}$.
\end{enumerate}\end{lemma}

\Proof We have $2v_1=-1$, $2v_2=1$, and $2b=6kp$, so $u_3=2b/(v\tau(v))=6kp$.
The prime $2$ being inert, by Lemma \ref{lemin2b}, if $4\nmid k$ the equation
is locally soluble at $2$ if and only if $\rho$ is a cube in $\F_4^*$, which
is not the case since the only cube is $1$. On the other hand, if $4\mid k$
Lemma \ref{lemin2b2} tells us that the equation is locally soluble at $2$.
Let us now look at the prime $p$. If $p\equiv1\pmod{3}$ then $p$ is split,
so by Corollary \ref{corsplit} the equation is locally soluble at $p$ if
and only if $(v_1+v_2d_p)X^3+(v_1-v_2d_p)Y^3+6kpZ^3=0$ is, and since the
$p$-adic valuations are $(0,0,\ge1)$, by Lemma \ref{lemu1u2u3} this is true if
and only if $(v_1+v_2d_p)/(v_1-v_2d_p)\equiv\rho^2/\rho=\rho\pmod{p}$ is
a cube in $\F_p^*$, hence if and only if $\rho^{(p-1)/3}\equiv1\pmod{p}$,
which is the case if and only if $p\equiv1\pmod{9}$. If $p\equiv2\pmod{3}$
then $p$ is inert, so by Lemma \ref{lemin2b} the equation is locally soluble
at $p$ if and only if $\rho$ is a cube in $\F_{p^2}^*$, hence if and only
if $p^2\equiv1\pmod{9}$, in other words $p\equiv-1\pmod{9}$ since we assume
$p\equiv2\pmod{3}$. It follows that the local condition at $p$ is
$p\equiv\pm1\pmod{9}$. Finally, let us look at the prime $3$. Since $2a=0$
we use Lemma \ref{lem3ram1} (2), which tells us that the equation is
locally soluble at $3$ if and only if $6kp\equiv\pm3\pmod{27}$, or
equivalently $kp\equiv\pm4\pmod{9}$, proving (1) since $p\equiv\pm1\pmod{9}$,
and (2) follows immediately.\fp

\smallskip

Next, we assume that $p\equiv1\pmod{3}$. In this case we can write
$p=\pi\tau(\pi)$ with $\pi=(w_1+w_2\sqrt{-3})/2$ in $12$ different
ways, and it is well-known and easy that up to sign and exchange of
$\pi$ and $\tau(\pi)$ there is exactly one such decomposition with $3\mid w_2$.

\begin{lemma} Let $p\equiv1\pmod{3}$, let $\pi=(w_1+w_2\sqrt{-3})/2$
be such that $\pi\tau(\pi)=p$, and let $k$ be such that $8\nmid k$,
$3\nmid k$, and $p\nmid k$. The equation corresponding to $\pi$ is locally 
soluble at the primes $2$, $3$, and $p$ if and only if either $4\mid k$ or 
$4\nmid k$ and $2\mid w_2\in\Z$, if $(w_1/(2k))^{(p-1)/3}\equiv1\pmod{p}$, and
if either $3\mid w_2$, or $3\nmid w_2$ and $p\equiv k^2+3\pmod{9}$.\end{lemma}

\Proof We have $2v_1=w_1$, $2v_2=w_2$, and $u_3=2b/(v\tau(v))=6k$.
The prime $2$ being inert, as above Lemmas \ref{lemin2b} and \ref{lemin2b2}
tell us that the equation is locally soluble at $2$ if and only if 
either $4\mid k$, or $4\nmid k$ and $(w_1+w_2\sqrt{-3})/2=1$ in $\F_4^*$,
which is equivalent to $2\mid w_2$. 
By Corollary \ref{corsplit} the equation is locally soluble at $p$ if
and only if $(v_1+v_2d_p)X^3+(v_1-v_2d_p)Y^3+6kd_pZ^3=0$ is. Since 
$v_1^2-v_2^2d_p^2=p$, we may assume for instance that $d_p$ is chosen so that
$v_p(v_1-v_2d_p)=1$, so in particular $v_2d_p\equiv v_1\pmod{p}$. The $p$-adic
valuations of the coefficients are thus $(0,1,0)$, so by Lemma \ref{lemu1u2u3}
local solubility is equivalent to $(v_1+v_2d_p)/(6kd_p)$ being a cube in 
$\F_p^*$, and since $v_1\equiv v_2d_p\pmod{p}$, this means that 
$v_1/(3kd_p)=w_1/(6kd_p)$ is a cube in $\F_p^*$. This is equivalent to
$(w_1/(6kd_p))^2=-w_1^2/(108k^2)$ being a cube, hence to $(w_1/(2k))^2$
being a cube, hence to $w_1/(2k)$ being a cube, leading to the given condition.
Finally, let us look at the prime $3$. Since $2a=0$ we use Lemma 
\ref{lem3ram1} (2), which tells us (since $3\nmid k$, so that $v_3(u_3)=1$) 
that the equation is locally soluble
at $3$ if and only if either $3\mid w_2$, or if $6k\equiv\pm 3w_1\pmod{27}$,
in other words $2k\equiv\pm w_1\pmod{9}$. However since $w_1^2+3w_2^2=4p$,
if $3\nmid w_2$ we have $w_1^2\equiv 4p-3\pmod{9}$, and since the condition
$2k\equiv\pm w_1\pmod{9}$ is equivalent to $w_1^2\equiv 4k^2\pmod{9}$ (since
$3\nmid w_1$) we obtain the equivalent condition $4p\equiv 4k^2+3\pmod{9}$,
or equivalently $p\equiv k^2+3\pmod{9}$, finishing the proof of the
lemma.\fp

\begin{corollary}\label{wh1} Let $p\ge5$ be prime, let $k=1$, $2$, or $4$,
and let $\wh{E}$ be the elliptic curve $y^2=x^3-27(kp)^2$.
\begin{enumerate}\item For $k=1$ or $k=2$, if either $p\equiv2\pmod{3}$, or 
if $p\equiv1\pmod{3}$, $p\not\equiv k^2+3\pmod{9}$, and 
$2^{(p-1)/3}\not\equiv1\pmod{p}$ then $\Ima(\wh{\al})$ is trivial.
\item For $k=4$, if $p\equiv2$ or $5\pmod{9}$ then $\Ima(\wh{\al})$ is trivial,
\item Otherwise, $|\Ima(\wh{\al})|=1$ or $3$ when $k=1$, $k=2$, or $k=4$
and $p\equiv8\pmod{9}$, and $|\Ima(\wh{\al})|=1$, $3$, or $9$ when $k=4$
and $p\equiv1\pmod{3}$.\end{enumerate}\end{corollary}

\Proof Since $2b=6kp$, $2$ is inert, and $3$ is ramified, with the notation
of Section \ref{secg3}, we must have $f_1=1$ if $p\equiv2\pmod{3}$ and
$f_1=1$ or $p$ if $p\equiv1\pmod{3}$. In the first case, the only possible
$v$ are $1$, $\rho$, and $\rho^2$, while in the second case we have in
addition the three possible $\pi$ (up to sign and conjugation) such that 
$\pi\tau(\pi)=p$. It follows that:

\begin{itemize}\item If $\rho\notin\Ima(\wh{\al})$ and none of the three 
possible $\pi$ is in $\Ima(\wh{\al})$ then $\Ima(\wh{\al})=\{1\}$, so
$|\Ima(\wh{\al})|=1$,
\item If $\rho\not\in\Ima(\wh{\al})$ and one of the three possible $\pi$ (so
necessarily exactly one) is in $\Ima(\wh{\al})$ then 
$\Ima(\wh{\al})=\{1,\pi,\tau(\pi)\}$, hence $|\Ima(\wh{\al})|=3$.
\item If $\rho\in\Ima(\wh{\al})$ and none of the three possible $\pi$
(in this case they are equivalent) is in $\Ima(\wh{\al})$ then
$\Ima(\wh{\al})=\{1,\rho,\rho^2\}$, hence $|\Ima(\wh{\al})|=3$.
\item If $\rho\in\Ima(\wh{\al})$ and one (hence all) of the three possible
$\pi$ are in $\Ima(\wh{\al})$ then
$\Ima(\wh{\al})=\{\rho^j,\rho^j\pi,\rho^j\tau(\pi),\ 0\le j\le2\}$, hence
$|\Ima(\wh{\al})|=9$.
\end{itemize}

Since in the two preceding lemmas we have studied local solubility of the 
equation in all these cases, we conclude by inspection.\fp

We can say a little more:

\begin{proposition}\label{propex2}
\begin{enumerate}\item Assume that $k=1$ and $p\equiv4\pmod9$,
or that $k=2$ and $p\equiv7\pmod9$, and write $p=m^2+3n^2$, where $m$ and $n$
are integers which are unique up to sign. The equation corresponding to
$\pi=m+n\sqrt{-3}$ (i.e., with $v_1=m$ and $v_2=n$) is ELS, and 
$|\Ima(\wh{\al})|=3$ if and only if it is globally soluble.
\item Assume that $k=4$ and $p\equiv 8\pmod{9}$. The equation corresponding
to $\rho$ is ELS, and $|\Ima(\wh{\al})|=3$ if and
only if it is globally soluble.
\item Assume that $k=4$ and $p\equiv 4$ or $7\pmod9$, and write
$4p=m^2+27n^2$, where $m$ and $n$ are unique up to sign. The equation 
corresponding to $\pi=(m+3n\sqrt{-3})/2$ (i.e., with $v_1=m/2$ and
$v_2=3n/2$) is ELS, and $|\Ima(\wh{\al})|=3$ if and 
only if it is globally soluble.
\end{enumerate}
\end{proposition}

\Proof Apply the same method as above. \fp

We will see below that it follows from BSD that these equations (of Proposition \ref{propex2}) should in
fact always be globally soluble.

\begin{corollary}\begin{enumerate}\item If $p\equiv5\pmod9$, or 
$p\equiv1$ or $7\pmod{9}$ and $2^{(p-1)/3}\not\equiv1\pmod{p}$, the elliptic 
curve $E_p$ with equation $y^2=x^3+p^2$ has rank $0$. If $p\equiv2\pmod9$,
it has rank $1$. Otherwise, if $p\equiv4$ or $8\pmod{9}$ it has rank $0$ or 
$1$, and if $p\equiv1$ or $7\pmod{9}$ it has rank $0$, $1$, or $2$.
\item If $p\equiv2\pmod{9}$, or $p\equiv1$ or $4\pmod{9}$ and 
$2^{(p-1)/3}\not\equiv1\pmod{p}$, the elliptic curve $E_{2p}$ with equation
$y^2=x^3+4p^2$ has rank $0$. Otherwise, if $p\equiv5$, $7$, or $8\pmod{9}$
it has rank $0$ or $1$, and if $p\equiv1$ or $4\pmod{9}$ it has rank $0$, $1$,
or $2$.
\item If $p\equiv2$ or $5\pmod{9}$ the elliptic curve $E_{4p}$ with equation
$y^2=x^3+16p^2$ has rank $0$. Otherwise, if $p\equiv4$, $7$, or $8\pmod{9}$ it
has rank $0$ or $1$, and if $p\equiv1\pmod{9}$ it has rank $0$, $1$, or $2$.
\end{enumerate}\end{corollary}

\Proof Clear, since $|\Ima(\al)||\Ima(\wh{\al})|=3^{r+1}$. (No need of BSD here.)\fp

This corollary allows us to determine the rank exactly for instance with $k=1$
for $p=61$, $79$, $113$, $131$, $149$, $151$, $163$, $293$, etc..., with $k=2$
for $p=29$, $83$, $137$, $139$, $173$, $181$, $199$, etc...,
and with $k=4$ for $p=41$, $59$, $101$, $131$, $137$, etc... for which 
{\tt mwrank}, at least in its basic version, is not able to determine the rank
using $2$-descent.

\smallskip

\noindent
{\bf Remarks.}\begin{enumerate}\item We can use ``the parity conjecture'' in this context, see for example \cite{Dok1} and \cite{Dok2},
in other words the analytic rank has the same parity as the
algebraic rank, then whenever in the above the rank is known to be equal to
$0$ or $1$ then it is always equal to $1$, while when the rank is known to
be equal to $0$, $1$, or $2$ then it is always equal to $0$ or $2$, and both
cases occur. This has been proved in certain cases: as already mentioned,
by Satg\'e for $k=1$ and $p\equiv2\pmod9$, and in an unpublished work Elkies 
has shown that for $k=4$ and $p\equiv4$ or $7\pmod{9}$ the rank is indeed 
equal to $1$.
\item The case $k=4$ corresponds to primes which are sums of two cubes,
so by \cite{RoZa} one knows that when $p\equiv1\pmod{9}$ the rank is equal
to $0$ or $2$, and exactly for which primes it is equal to $2$.
It is possible that either their method or Elkies' can be extended to the 
cases $k=1$ and $k=2$.
\item The result for $k=4$ can also be proved, less naturally, using
$2$-descent; see Theorem 6.4.17 of \cite{Cohen3}.
\end{enumerate}

\smallskip

When the cubics are ELS, we may of course try to look for a global solution
by search. A very efficient way of looking for rational points on a homogeneous
cubic has been described by N.~Elkies in \cite{Elkants}, see also an 
unpublished preprint of J.~Cremona
on the subject. It has been implemented by several people. Using a slightly
modified implementation due to M.~Watkins, we can for instance find that
a generator $P=(x,y)$ of the Mordell--Weil group of $y^2=x^3+p^2$ for 
$p=1759$, which has rank $1$, is given by

\begin{align*} x&=-\dfrac{242479559514608433100075350499874221113923535}{3063551062176562878606796987394973602467684}\\
y&=\dfrac{8643240396318605197724619647046515784779281219388876514209037894857}{5362134274928159502186511847328850266140274118035321166956948248}
\end{align*}

This generator is not found by {\tt mwrank} even at a high search limit. On the
other hand it could certainly be found using the Heegner point method.

For a more complicated example, for $p=9511$ the curve $y^2=x^3+p^2$ has
analytic rank $2$, so the Heegner point method is not applicable, and
{\tt mwrank} even at a high search limit finds only the one-dimensional
subspace of the (free part of the) Mordell--Weil group generated by
$P_1=(-210,9011)$. Using our implementation, we find that the full free
part has basis $(P_1,P_2)$ with $P_2=(x,y)$, where

\begin{align*} x&=\dfrac{32701984517186448621442294824950874787830128281}{456289760665179363242981599270033206574137600}\\
y&=\dfrac{92890043770264171014255964610503972850176417273682124237369198272789821}{9746778232027925565271633950191532413151456450450966045051557376000}\end{align*}

\medskip

In the following tables, we summarize what is proved (either using $3$-descent
as above, by Satg\'e, in Elkies' unpublished work, or 
Rodriguez-Villegas--Zagier's work),
what is a consequence of BSD, and what remains to be done. The tables are
coded as follows. In the first column we indicate the residue of $p$ modulo
$9$, and if relevant, in the second column we indicate the cubic character
$\lgc{2}{p}$ of $2$ modulo $p$, $1$ meaning that $2^{(p-1)/3}\equiv1\pmod{p}$,
and $\rho$, $\rho^2$ meaning of course $2^{(p-1)/3}\not\equiv1\pmod{p}$. In
the third, fourth, and fifth column we
give $|\Ima(\al)|$, $|\Ima(\wh{\al})|$, and the rank of the curve respectively,
and when two values are given, both occur. In the last column, we give
a pair of symbols (A,B), corresponding to $(|\Ima(\al)|,|\Ima(\wh{\al})|)$,
where P means proved, BSD means proved under BSD, S means Satg\'e, ELK means 
Elkies, RVZ means Rodriguez-Villegas--Zagier, and U means unknown.

\bigskip

\centerline{
\begin{tabular}{|l|l||l|l|l|l|}
\hline
$p\text{ mod 9}$ &$\lgc{2}{p}$  & $|\Ima(\al)|$ & $|\Ima(\wh{\al})|$ & \text{rank} & \text{ proved}\\ 
\hline\hline
$1$ & $1$ & $9\text{ or }3$ & $3\text{ or }1$ & $2\text{ or 0}$ & (U,U)\\
$1$ & $\rho$, $\rho^2$ & $3$ & $1$ & $0$ & (P,P)\\
$2$ & $-$ & $9$ & $1$ & $1$ & (S,P)\\
$4$ & $-$ & $3$ & $3$ & $1$ & (P,BSD)\\
$5$ & $-$ & $3$ & $1$ & $0$ & (P,P)\\
$7$ & $1$ & $9\text{ or 3}$ & $3\text{ or }1$ & $2\text{ or 0}$ & (U,U)\\
$7$ & $\rho$, $\rho^2$ & $3$ & $1$ & $0$ & (P,P)\\
$8$ & $-$ & $9$ & $1$ & $1$ & (BSD,P)\\
\hline
\end{tabular}}

\medskip

\centerline{Curves $y^2=x^3+p^2$, $p\ge5$}

\bigskip

\centerline{
\begin{tabular}{|l|l||l|l|l|l|}
\hline
$p\text{ mod 9}$ &$\lgc{2}{p}$  & $|\Ima(\al)|$ & $|\Ima(\wh{\al})|$ & \text{rank} & \text{ proved}\\ 
\hline\hline
$1$ & $1$ & $9\text{ or }3$ & $3\text{ or }1$ & $2\text{ or 0}$ & (U,U)\\
$1$ & $\rho$, $\rho^2$ & $3$ & $1$ & $0$ & (P,P)\\
$2$ & $-$ & $3$ & $1$ & $0$ & (P,P)\\
$4$ & $1$ & $9\text{ or 3}$ & $3\text{ or }1$ & $2\text{ or 0}$ & (U,U)\\
$4$ & $\rho$, $\rho^2$ & $3$ & $1$ & $0$ & (P,P)\\
$5$ & $-$ & $9$ & $1$ & $1$ & (BSD,P)\\
$7$ & $-$ & $3$ & $3$ & $1$ & (P,BSD)\\
$8$ & $-$ & $9$ & $1$ & $1$ & (BSD,P)\\
\hline
\end{tabular}}

\medskip

\centerline{Curves $y^2=x^3+4p^2$, $p\ge5$}

\bigskip

\centerline{
\begin{tabular}{|l|l||l|l|l|l|}
\hline
$p\text{ mod 9}$ &$\lgc{2}{p}$  & $|\Ima(\al)|$ & $|\Ima(\wh{\al})|$ & \text{rank} & \text{ proved}\\ 
\hline\hline
$1$ & $-$ & $3$ & $9\text{ or }1$ & $2\text{ or 0}$ & (P,RVZ)\\
$2$ & $-$ & $3$ & $1$ & $0$ & (P,P)\\
$4$ & $-$ & $3$ & $3$ & $1$ & (P,ELK)\\
$5$ & $-$ & $3$ & $1$ & $0$ & (P,P)\\
$7$ & $-$ & $3$ & $3$ & $1$ & (P,ELK)\\
$8$ & $-$ & $3$ & $3$ & $1$ & (P,BSD)\\
\hline
\end{tabular}}

\medskip

\centerline{Curves $y^2=x^3+16p^2$, $p\ge5$}

\bigskip

An immediate corollary of the above tables and of Corollary \ref{corex1} and
Proposition \ref{propex2} is the following:

\begin{corollary} \begin{enumerate}\item Assume BSD. If $p\equiv2$ or 
$8\pmod9$ there exist $x$ and $y$ in $\Q$ such that $p=x^3+2y^3$, and if 
$p\equiv5$ or $8\pmod9$ there exist $x$ and $y$ in $\Q$ such that $p=x^3+4y^3$.
\item Assume that either $k=1$ and $p\equiv4\pmod{9}$ or that $k=2$ and
$p\equiv7\pmod{9}$, and write $p=m^2+3n^2$. If BSD is true the equation
$$nX^3-3mY^3+3kZ^3+3mX^2Y-9nXY^2=0$$
is globally soluble.
\item Assume that $p\equiv8\pmod{9}$. If BSD is true the equation
$$X^3+3Y^3+24pZ^3-3X^2Y-9XY^2=0$$
is globally soluble.
\item Assume that $p\equiv4$ or $7\pmod9$, and write $4p=m^2+27n^2$. Without
any assumption the equation
$$nX^3-mY^3+8Z^3+mX^2Y-9nXY^2=0$$
is globally soluble.
\end{enumerate}
\end{corollary}

\Proof Clear, since these correspond respectively to $|\Ima(\al)|=9$ under
BSD, $|\Ima(\wh{\al})|=3$ under BSD for (2) and (3), and to 
$|\Ima(\wh{\al})|=3$ by Elkies's result.\fp

\subsection{The Curves $y^2=x^3+(kp)^2$ for $k=3$ or $9$}

Once again the restriction on $k$ is made so that no other prime apart from
$3$ divides it, and it is not necessary to consider higher powers of $3$.
Furthermore the primes $p=2$ and $p=3$ give rise to a finite number of curves
which can be treated individually (specifically, the rank is zero unless
$(k,p)=(3,2)$, in which case it has rank $1$, a Mordell--Weil generator being
$(-3,3)$, and the torsion is of order $3$ generated by $T=(0,kp)$, unless
$(k,p)=(9,3)$, in which case it has order $6$ generated by $(18,81)$). We
therefore assume that $p\ge5$.

As before, we first compute the image of $\al$. For this, we consider the 
cubic equations of Theorem \ref{thd1} (3), in other words 
$u_1X^3+u_2Y^3+u_3Z^3=0$, where $u_1u_2\mid 2kp$ and $u_3=2kp/(u_1u_2)$, where
we recall that $u_1u_2$ is squarefree. Up to exchange of $u_1$ and $u_2$, it 
is easy to check that the only of possibilities are $(1,1,2kp)$, $(1,2,kp)$,
$(1,3,2kp/3)$, $(1,6,kp/3)$, $(1,p,2k)$, $(1,2p,k)$, $(1,3p,2k/3)$, 
$(1,6p,k/3)$, $(2,3,kp/3)$, $(2,p,k)$, $(2,3p,k/3)$, $(3,p,2k/3)$,
$(3,2p,k/3)$, and $(6,p,k/3)$. The first one (corresponding to $u=1$) gives
an evidently soluble equation, corresponding to the unit element of the
elliptic curve.

\smallskip

Consider first $k=3$. We obtain the equations
$(1,2,3p)$, $(1,3,2p)$, $(1,6,p)$, $(1,p,6)$, $(1,2p,3)$, $(1,3p,2)$,
$(1,6p,1)$, $(2,3,p)$, $(2,p,3)$, $(2,3p,1)$, $(3,p,2)$, $(3,2p,1)$, and
$(6,p,1)$. The equation $(1,6p,1)$ (corresponding to $u=6p$ and $u=(6p)^2$)
corresponds to the two $3$-torsion points of the curve. Apart from that, up
to permutation of the $u_i$ we have to study solubility for
$(1,2,3p)$ (corresponding to $u=2$, $4$, $3p$, $9p^2$, $12p$, and $18p^2$),
$(1,3,2p)$ (corresponding to $u=3$, $9$, $2p$, $4p^2$, $18p$, and $12p^2$),
$(1,6,p)$ (corresponding to $u=6$, $36$, $p$, $p^2$, $36p$, and $6p^2$),
$(2,3,p)$ (corresponding to $u=12$, $18$, $4p$, $2p^2$, $9p$, and $3p^2$).

\begin{lemma} We have two cases :
\begin{enumerate}\item If $p\equiv2\pmod{3}$ all the above
cubics are ELS, giving a total of $27$ ELS cubics.
\item If $p\equiv1\pmod{3}$, then either both $2$ and $3$ are cubes in
$\F_p^*$, in which case once again all the above cubics are ELS
for a total of $27$, or either $2$ or $3$ or both are non cubes,
in which case only $9$ equations are ELS.
\end{enumerate}
\end{lemma}

\Proof Using the lemmas that we have proved it is immediate to show that
all the equations are soluble at $2$ and $3$. The only problem is at $p$,
and by Lemma \ref{lemu1u2u3} the equations are also soluble at $p$ if and
only if the classes of $2$, $3$, $6$, and $3/2$ respectively are
cubes in $\F_p^*$, and if $p\equiv2\pmod{3}$ this is trivially true.

Assume now that $p\equiv1\pmod{3}$. We consider four cases, according to the
cubic residue character of $2$ and $3$ modulo $p$.

\begin{enumerate}\item If $2$ and $3$ are non cubes in $\F_p^*$, then
either their product or their quotient is a cube, so we deduce that either
$(1,6,p)$ or $(2,3,p)$ (but not both) is locally soluble, giving a total
of $6+3=9$ possible values of $u$.
\item If $2$ or $3$ is a cube in $\F_p^*$ but not both, then $6$ and $3/2$
cannot be cubes, so we deduce that either $(1,2,3p)$ or $(1,3,2p)$ is locally
soluble, giving again $9$ possible values of $u$.
\item If both $2$ and $3$ are cubes in $\F_p^*$ then all the equations
are locally soluble, giving a total of $24+3=27$ possible values of $u$,
proving the lemma.\fp
\end{enumerate}

\begin{corollary}\begin{enumerate}\item If $p\equiv2\pmod{3}$ then
$|\Ima(\al)|=3$, $9$, or $27$ and assuming BSD we have $|\Ima(\al)|=3$ or
$27$ and the following are equivalent:
\begin{enumerate}\item $|\Ima(\al)|=27$,
\item $\rk(E)=2$,
\item There exist $x$ and $y$ in $\Q$ such that $p=x^3+6y^3$,
\item There exist $x$ and $y$ in $\Q$ such that $p=2x^3+3y^3$,
\item There exist $x$ and $y$ in $\Q$ such that $p=4x^3+12y^3$,
\item There exist $x$ and $y$ in $\Q$ such that $p=9x^3+18y^3$.
\end{enumerate}
\item If $p\equiv1\pmod{3}$ and $2$ and $3$ are not both cubes
in $\F_p^*$ then $|\Ima(\al)|=3$ or $9$, and assuming BSD we have
$|\Ima(\al)|=9$ and $\rk(E)=1$, and furthermore for $(a,b)=(1,6)$, $(2,3)$, 
$(4,12)$, or $(9,18)$, there exist $x$ and $y$ in $\Q$ such that $p=ax^3+by^3$
if and only if $b/a$ is a cube in $\F_p^*$.
\item If $p\equiv1\pmod{3}$ and $2$ and $3$ are both cubes in $\F_p^*$ then
either $C_{\pi_0}$ is globally soluble, in which case the six conditions
of (1) are again equivalent except that the second must be replaced by
$\rk(E)=3$, or $C_{\pi_0}$ is not globally soluble, in which case 
$|\Ima(\al)|=9$ and $\rk(E)=1$, and exactly one of the equations
$p=x^3+6y^3$, $p=2x^3+3y^3$, $p=4x^3+12y^3$, or $p=9x^3+18y^3$ has a solution
with $x$ and $y$ in $\Q$.
\end{enumerate}\end{corollary}
%it follows from Lemma \ref{lemwal}
\Proof The assertions independent of BSD follow immediately from the above
lemma. On the other hand, we have that if either
$p\equiv2\pmod3$ or $p\equiv1\pmod{3}$ and $2$ and $3$ are not both cubes in
$\F_p^*$ we have $|\Ima(\wh{\al})|=1$, so that $3^{\rk(E)+1}=|\Ima(\al)|$.
If $p\equiv2\pmod{3}$ the root number of $E$ is equal to $+1$, so assuming
BSD the rank of $E$ is even, so we have $|\Ima(\al)|=3$ or $27$, and the
result clearly follows in this case. If $p\equiv1\pmod{3}$ the root number of
$E$ is equal to $-1$, so assuming BSD the rank of $E$ is odd, so we have
$|\Ima(\al)|=9$ when $2$ and $3$ are not both cubes, so the result also 
follows. Finally, if $p\equiv1\pmod{3}$ and $2$ and $3$ are both cubes, then
$|\Ima(\wh{\al})|=1$ or $3$, and it is equal to $3$ if and only if $C_{\pi_0}$
(which is ELS) is globally soluble. Assuming BSD the rank of $E$ is again odd,
so we have two cases:

$\bullet$ If $C_{\pi_0}$ is not globally soluble we again have 
$|\Ima(\wh{\al})|=1$, so we must have $|\Ima(\al)|=9$ and $\rk(E)=1$.

$\bullet$ If $C_{\pi_0}$ is globally soluble we have $|\Ima(\wh{\al})|=3$,
so we must have either $|\Ima(\al)|=3$ and $\rk(E)=1$, or $|\Ima(\al)|=27$ and
$\rk(E)=3$.\fp

Note that although $C_{\pi_0}$ is ELS, it is not always globally soluble:
the smallest $p$ for which it is not is $p=3889$, for which the $2$-rank based
{\tt mwrank} program of Cremona tells us that the rank is equal to $1$,
and on the other hand $(X,Y,Z)=(91,-211,19)$ is a solution of the $(2,3,p)$
cubic so $|\Ima(\al)|\ge9$, hence by the above if $C_{\pi_0}$ was globally
soluble we would have $\rk(E)=3$, which is not the case.

\smallskip

Consider now $k=9$. We obtain the equations $(1,2,9p)$, $(1,3,6p)$, $(1,6,3p)$,
$(1,p,18)$, $(1,2p,9)$, $(1,3p,6)$, $(1,6p,3)$, $(2,3,3p)$, $(2,p,9)$,
$(2,3p,3)$, $(3,p,6)$, $(3,2p,3)$, and $(6,p,3)$. The equation $(3,2p,3)$
(corresponding to $u=18p$ and $u=12p^2$) corresponds to the two $3$-torsion
points of the curve. Apart from that, up to permutation of the $u_i$ we have
to study solubility for $(1,2,9p)$ (corresponding to $u=2$ and $4$),
$(1,3,6p)$ (corresponding to $u=3$, $9$, $6p$, and $36p^2$),
$(1,6,3p)$ (corresponding to $u=6$, $36$, $3p$, $9p^2$),
$(1,9,2p)$ (corresponding to $u=2p$ and $4p^2$),
$(1,18,p)$ (corresponding to $u=p$ and $p^2$),
$(2,3,3p)$ (corresponding to $u=12$, $18$, $12p$, and $18p^2$),
$(2,9,p)$ (corresponding to $u=4p$ and $2p^2$), and
$(3,6,p)$ (corresponding to $u=9p$, $3p^2$, $36p$, and $6p^2$).

Once again, the equations are all locally soluble at $2$. However, they are
not all locally soluble at $3$ (in fact $(3,6,p)$ is never locally soluble
at $3$), and using once again the local solubility results that we have
proved, we obtain the following lemma:

\begin{lemma}\begin{enumerate}\item If $p\equiv2\pmod{3}$, exactly two of the
$7$ above equations are ELS, giving always a total of $9$ values of $u$.
\item If $p\equiv1\pmod3$ and $3/2$, $3$, or $6$ are cubes respectively for
$p\equiv1$, $4$, or $7\pmod9$, once again exactly two of the $7$ above
equations are ELS, giving always a total of $9$ values of $u$. Otherwise,
none are ELS, giving a total of $3$ values of $u$ (corresponding to the
$3$-torsion points).\end{enumerate}\end{lemma}

\Proof The seventh equation is not locally soluble at $3$, and the six others
are ELS if and only if $p\equiv\pm4\pmod{9}$ and $3$ is a cube for $(1,3,6p)$
and $(1,9,2p)$, $p\equiv\pm2\pmod{9}$ and $6$ is a cube for $(1,6,3p)$ and
$(2,9,p)$, or $p\equiv\pm1\pmod{9}$ and $3/2$ is a cube for $(1,18,p)$ and
$(2,3,3p)$, and the result follows.\fp

\bigskip

\nocite{*}

\bibliographystyle{amsalpha}

%\bibliographystyle{alpha}
%\bibliography{ref3descente}

\begin{thebibliography}{xx}
 \bibitem{Cas} \textsc{Cassels, J. W. S.}, 
\textit{Arithmetic on curves of genus {$1$}. {I}. {O}n a conjecture of
              {S}elmer}.
J. Reine Angew. Math., {\bf202} (1959), 52--99.

\bibitem{Cohen1} \textsc{Cohen, H.}, 
\textit{A course in computational algebraic number theory}.
Graduate Texts in Mathematics, Springer (New York), {\bf138} (1993).

\bibitem{Cohen2} \textsc{Cohen, H.}, 
\textit{Advanced topics in computational number theory}.
Graduate Texts in Mathematics, Springer (New York), {\bf193} (2000).

\bibitem{Cohen3} \textsc{Cohen, H.}, 
\textit{Number theory. {V}ol. {I}. {T}ools and {D}iophantine
              equations}.
Graduate Texts in Mathematics, Springer (New York), {\bf239} (2007).

\bibitem{Cohen4} \textsc{Cohen, H.}, 
\textit{Number theory. {V}ol. {II}. {A}nalytic and modern tools}.
Graduate Texts in Mathematics, Springer (New York), {\bf239} (2007).

\bibitem{CreFiOSiSto1} \textsc{Cremona, J. E. and Fisher, T. A. and O'Neil, C. and Simon, D. and Stoll, M.}, 
\textit{Explicit n-descent on elliptic curves. I. Algebra}.
J. reine angew. Math., {\bf615} (2008), 121-155.

\bibitem{CreFiOSiSto2} \textsc{Cremona, J. E. and Fisher, T. A. and O'Neil, C. and Simon, D. and Stoll, M.}, 
\textit{Explicit n-descent on elliptic curves. II. Geometry}.
J. reine angew. Math., {\bf} (accepted in 2008).

\bibitem{DeLo} \textsc{DeLong, M.}, 
\textit{A formula for the {S}elmer group of a rational three-isogeny}.
Acta Arithmetica, {\bf105} (2002), 119--131.

\bibitem{DjaShaSma} \textsc{Djabri, Z. and Schaefer, E. F. and Smart, N. P.}, 
\textit{Computing the {$p$}-{S}elmer group of an elliptic curve}.
Trans. Amer. Math. Soc., {\bf352} (2000), 5583--5597.

\bibitem{Dok1} \textsc{Dokchitser, T. and Dokchitser V.}, 
\textit{Parity of ranks for elliptic curves with a cyclic isogeny}.
J. Number Theory, {\bf128} (2008), 662-679.

\bibitem{Dok2} \textsc{Dokchitser, T. and Dokchitser V.}, 
\textit{On the Birch-Swinnerton-Dyer quotients modulo squares}.
to appear in Annals of Math.

\bibitem{Elkants} \textsc{Elkies, N. D.}, 
\textit{Rational points near curves and small nonzero {$\vert x\sp
              3-y\sp 2\vert $} via lattice reduction}.
Algorithmic number theory ({L}eiden), Lecture Notes in Comput. Sci. {\bf1838} (2000), 33--63.

\bibitem{Fish} \textsc{Fisher, T.}, 
\textit{Finding rational points on elliptic curves using 6-descent and
              12-descent}.
J. Algebra, {\bf320} (2008), 853--884.

\bibitem{Lem} \textsc{Lemmermeyer, F.}, 
\textit{Reciprocity laws}.
Springer Monographs in Mathematics, From Euler to Eisenstein {\bf1838} (2000).

\bibitem{RoZa} \textsc{Rodr{\'{\i}}guez Villegas, F. and Zagier, D.}, 
\textit{Which primes are sums of two cubes?}.
Number theory (Halifax, NS), CMS Conf. Proc. {\bf15} (1994), 295--306.

\bibitem{Satg} \textsc{Satg{\'e}, P.}, 
\textit{Groupes de {S}elmer et corps cubiques}.
J. Number Theory {\bf23} (1986), 294--317.

\bibitem{Satg2} \textsc{Satg{\'e}, P.}, 
\textit{Un analogue du calcul de {H}eegner}.
Invent. Math. {\bf87} (1987), 425--439.

\bibitem{SchaStoll} \textsc{Schaefer, E. F. and Stoll, M.}, 
\textit{How to do a {$p$}-descent on an elliptic curve}.
Trans. Amer. Math. Soc. {\bf356} (2004), 1209--1231 (electronic).

\bibitem{Sil1} \textsc{Silverman, J. H.}, 
\textit{The arithmetic of elliptic curves}.
Graduate Texts in Mathematics, Springer-Verlag (New York) {\bf106} (1992).

\bibitem{Top} \textsc{Top, J.}, 
\textit{Descent by {$3$}-isogeny and {$3$}-rank of quadratic fields}.
Advances in number theory (Kingston, ON), Oxford Univ. Press {\bf106} (1991).


\end{thebibliography}

\end{document}